%% file: paper.tex
\documentclass[a4paper, 
11pt
]{article}
\usepackage{xcolor}
\usepackage{mathtools}
\usepackage{symbols}
\usepackage{amsfonts}
\usepackage{xspace}
\usepackage{empheq}
\usepackage{cancel}
\usepackage{url}  
\usepackage{rotating} 
\usepackage{placeins}

\usepackage{ifpdf}

\ifpdf
\usepackage{pdfsync}
\fi
\usepackage{todonotes}

\newcommand{\phins}{{{\frac{\partial\phi}{\partial n}}}}
\newcommand{\deal}{\texttt{deal.II}\xspace}

\begin{document}

\title{A unified steady and unsteady formulation for hydrodynamic potential flow simulations with
       fully nonlinear free surface boundary conditions}

\author{Andrea~Mola\footnote{Scuola IMT Alti Studi Lucca, Piazza S. Ponziano, 6 - 55100 Lucca, Italy, andrea.mola@imtlucca.it},
        Nicola~Giuliani\footnote{previously Scuola Internazionale Sueriore di Studi Avanzati (SISSA), now Applied Materials Inc., Santa Clara, CA 95054-3299 USA, nicola\_giuliani@amat.com},
        \'Oscar~Crego\footnote{Departamento de Matemática Aplicada, Universidade de Santiago de Compostela, Santiago de Compostela E-15782, Spain. oscar.crego@usc.es},
       Gianluigi~Rozza\footnote{Scuola Internazionale Sueriore di Studi Avanzati (SISSA) via Bonomea, 265 - 34136 Trieste, Italy, grozza@sissa.it}}

\maketitle

\begin{abstract}
This work discusses the correct modeling of the fully nonlinear free surface
boundary conditions to be prescribed in water waves flow simulations
based on potential flow theory. The main goal of such a discussion is
that of identifying a mathematical formulation and a numerical treatment
that can be used both to carry out transient simulations, and to
compute steady solutions --- for any flow admitting them. In the literature
on numerical towing tank in fact, steady and unsteady fully nonlinear
potential flow solvers are characterized by different mathematical formulations.
The kinematic and dynamic fully nonlinear free surface boundary conditions are
discussed, and in particular it is proven that the kinematic free surface
boundary condition, written in semi-Lagrangian form, { can be manipulated to derive
an alternative} non penetration boundary condition by all means identical
to the one used on the surface of floating bodies or on the basin bottom.
The simplified mathematical problem obtained is discretized over space and time via
Boundary Element Method (BEM) and Implicit Backward Difference Formula (BDF)
scheme, respectively. The results confirm that the solver implemented is
able to solve steady potential flow problems just by eliminating
null time derivatives in the unsteady formulation. Numerical results obtained
confirm that the solver implemented is able to accurately
reproduce results of classical steady flow solvers available in the
literature.\\

\noindent
The final version of the present paper has been accepted for publication on Applied Mathematical Modelling.
\end{abstract}

\tableofcontents
\section{Introduction and literature review}\label{sec:intro}

The progress witnessed in the last decades has established
computational tools for fluid dynamic performance prediction
as a reliable instrument available to boat and ship designers, 
and a valid alternative to the experimental approach.
Along with a steady increase in computational power and resources,
such progress has to be ascribed to the constant improvement of
mathematical models and numerical algorithms. Among the many methods
developed in the effort to obtain fast and yet accurate hydrodynamic
simulations, potential flow models complemented by fully nonlinear
free surface boundary conditions have enjoyed considerable success
in the naval architecture community. In fact both the incompressible
fluid and irrotational flow assumptions upon which the potential
flow theory is based appear quite reasonable for slender hulls
advancing at moderate cruise speeds. In addition, compared to their
linearized free surface boundary condition counterparts, fully nonlinear
potential models enjoy superior accuracy, which makes them able to
predict displacement hulls resistance with errors as low as 2\% and
water elevations within experimental uncertainty \cite{molaIsope2016}.
On the other hand, compared to more general models based on Navier--Stokes
equations, such as RANS or LES, they clearly lead to smaller discretized
problems and to faster computations.

One of the most important traits of the flow past a ship hull is given by the
presence of two fluids --- air and water --- around it, separated by a sharp interface, or \emph{free surface}.
Correct modeling of the free surface wave pattern surrounding the hull is paramount for
accurate estimation of the energy dispersed by the ship to the surrounding waves,
and ultimately to a good estimate of the hydrodynamic forces. This is in essence
the reason why fully nonlinear potential models are able to provide accurate fluid dynamic forces
predictions despite the significant simplifying assumptions --- irrotational flow,
inviscid fluid and simply connected domain --- upon which they are based.
In the framework of potential flow models with fully nonlinear free surface conditions,
the governing Laplace equation for the velocity potential is only solved in the
portion of space surrounding the hull and occupied by water. Thus, the position and shape
assumed at each time instant by the free surface is an additional unknown of the
resulting mathematical problem. Along with the boundary condition on the velocity potential, an
additional boundary condition must be added to compute the evolution of the free surface position.
 
The most common approach for  unsteady free surface potential
flow simulations (see, e.g., Grilli et al. \cite{grilli2001}) is the mixed Eulerian–Lagrangian approach (MEL) originally
introduced by Longuet-Higgins and Cokelet \cite{longuet-higginsCokeletMEL76}.
In such a framework, at each time step
a Laplace boundary value problem for the fluid velocity potential is solved
in the Eulerian step making use of a Dirichlet free surface boundary condition. 
The resulting fluid velocity field is then introduced into suitable kinematic
and dynamic free surface conditions to compute, in a Lagrangian step, the time
evolution of the free surface position and potential to be used at the next
time step. In the staggered time integration approach characterizing MEL formulation,
the fluid dynamic grid nodes follow in a Lagrangian fashion the fluid particles
on the free surface. In typical environmental applications this is normally
not a thing of concern, as water waves are associated with small mass transport
as well as small average particle velocities. However, ship hydrodynamics simulations
are usually carried out in the frame of reference of the moving hull as it advances
through the water. In such a frame, the presence of a stream flowing past the
ship means that the application MEL leads to undesired downstream drift of the mesh
nodes, which would move the  numerical domain away from the region of
interest --- the hull surroundings --- as the simulation is carried out. For such a
reason, the use of MEL in ship hydrdodynamics applications requires periodic
regridding or complex grid treatments, which increases the computational cost
and the implementation complexity of the
algorithm (see, e.g., Kjellberg et al. \cite{kjellbergEtAl2011}, and
Kjellberg \cite{kjellbergPhD2013}). In addition, it must be remarked that the
hull frame of reference is also the only one in which a steady state flow is
obtained. As well known, a flow is steady when the fluid properties at each
point in the domain do not change over time. Indeed, the fact that a flow can have
a steady description can depend on the chosen frame of reference. In the
case at hand, the flow past a moving hull might admit a steady description
only in the frame of reference of the moving hull, in which the hull itself 
(and the boundary associated with it) is stationary. In any frame in which the
hull position depends on time, a stationary flow solution cannot instead be identified.
Clearly, if a steady flow exists, steady solvers designed to take advantage of
the absence of time dependence can find the flow solution solving a single non
linear problem. This represents a clear computational advantage compared to the multiple
problems each associated with a time step of time dependent problems, which have to be solved
until a steady state solution is approximately reached.

For such reasons, researchers in naval architecture community have always shown great
interest in developing steady solvers for the potential flow model with fully nonlinear
free surface boundary conditions. However, given the problems in its application
in the hull reference frame, MEL cannot be used in its original formulation to solve
these problems. More specifically, even in presence of steady state flow, the free surface
nodes position computed by MEL is not steady but follows the water motion downstream.
Simply dropping time derivatives in the spatially discretized numerical problem resulting from 
MEL does not lead to numerical problems with steady solution. Despite this difficulties
several researchers obtained converging algorithms for the solution of steady fully nonlinear
free surface potential flows. A number of different potential flow models are available in the literature
for solving steady nonlinear free surface flows past a ship hull.
Among others, we mention the work of Raven \cite{ravenPhD1998}, which resulted in the
implementation of the commercial software RAPID, and of Janson \cite{JansonPhD1997,janson2000},
which led to the commercial code SHIPFLOW, and of Scullen \cite{Scullen}. Over the years, these
steady state solvers have established themselves as fast and reliable
tools for the early design stages, in which they can provide not only
ship wave resistance estimates, but also pressure distributions,
free surface elevation and velocity fields surrounding the hull.

Such algorithms are quite different one from each other. In Raven's work \cite{ravenPhD1998}
the main idea is combining together the two free surface boundary conditions used in MEL, so as
to obtain a single free surface condition in which time derivatives are dropped. Scullen \cite{Scullen}
instead uses a non penetration boundary condition on the idle free surface, and then updates its position
based on a dynamic condition where again Eulerian time derivatives are dropped. In all cases, an iterative
scheme based on the previous steps is used to update the free surface position and potential flow solution
until convergence. Despite its remarkable
effectiveness and accuracy in obtaining solutions for steady flows in a small number of iterations, this
kind of approach cannot be effectively used in presence of transient flows. As noted by Raven in the
introduction of his PhD dissertation (\cite{ravenPhD1998},page 67), a unified approach for steady and transient
fully nonlinear potential flow simulations, was missing at the time. And, to the best of the authors knowledge,
it is missing to this day. This fact represents a clear anomaly with respect to other fluid dynamic or more in
general continuum mechanics models. In such dynamical systems in fact, the steady solution is typically
sought simply through elimination of the time derivatives from the unsteady governing equations.

There is however a different approach that can be used to make the MEL approach more suitable
to moving reference frames in which non negligible fluid stream velocities are
observed. As the experience in the Finite Element Methods (FEM) community suggests,
the Arbitrary Lagrangian--Eulerian (or ALE, \cite{doneaetAlAle2004}) formulation
is an effective approach in dealing with moving boundaries in presence when significant
transport velocities. First attempts to employ a similar methodology  
have been carried out by Beck \cite{beck1994}, which developed a set of fully nonlinear
free surface boundary conditions written in semi-Lagrangian form.

 The time
derivatives appearing in Beck boundary conditions equations are neither computed
on fixed spatial points as in the Eulerian formulation, nor on fixed fluid
particles as in the Lagrangian formulation. Instead, they are computed on fixed
grid nodes --- or free surface markers ---, which move with a user prescribed
velocity field. Such a formulation, by all means similar to  ALE,
allows for the resolution of the nonlinear free surface problem using a
time advancing strategy identical to the one used in MEL,
and in principle it does not require regridding. However, a saw-
tooth instability arising from dominant transport terms appearing in the newly developed
boundary conditions, make semi-Lagrangian unstable whenever the grid
and fluid velocity difference is non negligible. For such a reason, the
methodology combining semi-Lagrangian and MEL could only be used in \cite{scorpioPhD1997}
by setting a grid stream velocity equal to the fluid one, which only
mitigates the regridding problem. 
It was not until a decade ago (see \cite{waveBem,Giuliani2015}) that a proper stabilization
mechanism was introduced to allow for stable simulations based on semi-Lagrangian
boundary conditions, at all non planing boat speeds and with no remeshing required on unstructured and
adaptively refined grids. Moreover, \cite{waveBem} presents
a novel time implicit time advancing scheme as opposed to the explicit staggered approach
typically used with MEL.
The stabilized semi-Lagrangian free surface model presented in \cite{waveBem}, features
no average grid stream velocity and can in principle be used to obtain a solver that is
switched from unsteady to steady by only removing the time derivatives
from the governing equations. However,  at the numerical level the nonlinear problem resulting from such operation is not
able to converge to a steady state solution. For such reason, in this work, we present all the
modifications to the stabilized semi-Lagrangian free surface model
presented in \cite{waveBem}, so that it can be successfully used both as a steady and unsteady solver.
We report a table summarizing the features of the most common numerical approaches in Table~\ref{tab:techniques_comparison}.

As will be discussed, further modifications are made to the transient free surface model, to make it
also compatible with steady state solution. In particular, the semi-Lagrangian kinematic free surface
boundary condition is replaced by a non homogeneous Neumann boundary condition, written in ALE formulation.
Such a condition is substantially a non penetration constraint for the fluid on the moving boundary,
and is identical to the conditions imposed on other non penetration regions such as the hull or bottom
boundaries. The presence in the system of the dynamic semi-Lagrangian free surface boundary condition
allows for the simultaneous computation of the free surface grid velocity, which is the additional
unknown of the problem. A proof will be offered that the non penetration Neumann boundary condition
can be derived the semi-Lagrangian kinematic free surface condition on the free surface. Several
numerical experiments will then show that the transient solver can be used to obtain steady solution
just ``turning off'' all time derivatives in the governing equations. In addition, to confirm that the
present approach recovers the results of classic steady solvers, the steady state results will be compared
to results on the same test cases obtained by Scullen \cite{Scullen}.

The content of this paper is organized as follows. Section 2 introduces the equations of the
model for the free surface flow, based on the potential flow theory. Details of the free surface
modeling will be also presented in such section. Section 3 describes the numerical discretization of
the problem based on a combined Boundary Element Method (BEM) and Finite Element Method (FEM) approach, with
implicit Backward Difference Formula (BDF) time advancing scheme. Section 4 provides a description
of the numerical test cases considered, and of the results obtained. Finally, Section 5 reports some
brief conclusion remarks, and possible follow up investigations.

\section{Fluid dynamic model based on potential flow theory}

As mentioned, the main goal of this work is that of identifying a unified
mathematical model to carry out both steady and unsteady
simulations of the potential flow past a body advancing at constant speed
in calm water. We point out that by \emph{steady flow} we indicate a flow in
which all Eulerian derivatives of the unknown pressure field and velocity
field are null, namely

\begin{equation}
\dsfrac{\Pa p(\xb,t)}{\Pa t} = 0 \qquad \dsfrac{\Pa \ub(\xb,t)}{\Pa t} = 0.
\end{equation}
In such conditions we have that the flow fields $p(\xb,t)=p(\xb)$ and $\ub(\xb,t)=\ub(\xb)$
are only dependent on the point position vector $\xb=x\eb_x+y\eb_y+z\eb_z$, 

where $\eb_x$, $\eb_y$ and $\eb_z$, are the unit vectors associated with the axes of the chosen reference frame.

As stressed before the only frame of reference in which a steady state solution --- if any exists --- can be observed is a body attached one. In fact,
it is the only frame of reference in which the body surface, which is a boundary of the fluid domain, has a stationary position.
For such a reason, in this work we will describe the flow field in a frame of reference attached to the boat hull as it advances in the water. 
Of course, there are situations in which even in a body attached reference frame a steady
state is not possible. This is for instance the case of  the flow past a ship advancing
through waves in an unsteady wave field. In all those situations, resorting to an unsteady formulation will be mandatory.
In the present section we will start discussing the latter, more general case, and
successively consider its steady variant. 

\subsection{Governing equations}

\begin{center}
\begin{figure}
\centerline{
  \ifpdf
  \resizebox{0.7\textwidth}{!}{
    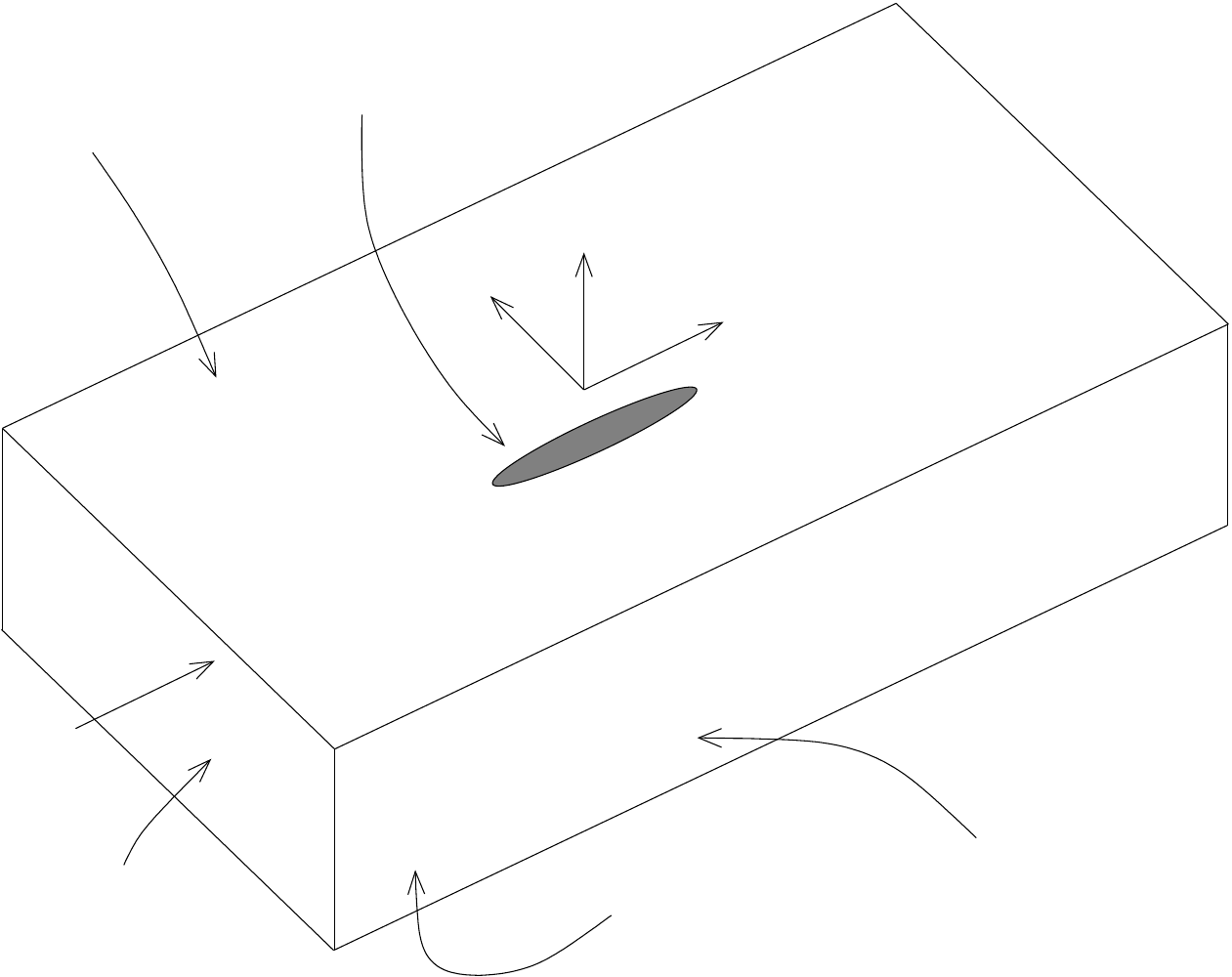
  }
  \else
  \resizebox{0.7\textwidth}{!}{
    \input{./figures/dominio_barca.pstex_t}
  }
  \fi
}
\caption{A sketch of the computational domain $\Omega(t)$. The diagram also shows the
         location of the body, free surface, bottom, inflow and far field boundary regions
         $\Gamma^h,\Gamma^{fs},\Gamma^b,\Gamma^i$ and $\Gamma^\infty$, respectively. 
         The far field stream velocity $\nablab\phi_\infty$ also indicates that $\Omega(t)$ is located in
         a hull attached frame, and follows the body in its motions.\label{fig:domain_and_regions}}
\end{figure}
\end{center}

The flow domain $\Omega(t),\ \forall t\geq 0$ is the simply connected and time dependent three dimensional
region occupied by water, surrounding and following the body mean flow velocity.
In such a region --- depicted in Figure \ref{fig:domain_and_regions}, along with its 
boundaries --- assuming irrotational flow, the velocity field admits a scalar potential
$\Phi(\xb,t)$, namely

\begin{equation}\label{eq:pert_pot}
\ub(\xb,t) = \nablab \Phi(\xb,t) = \nablab\phi_\infty(t) + \nablab\phi(\xb,t) 
\qquad \forall \xb=x\eb_x+y\eb_y+z\eb_z\ \in\ \Omega(t) \subseteq \mathbb{R}^3
\end{equation}
where $\phi(\xb,t)$ is the perturbation potential and  $\phi_\infty(\xb,t)$ is the
asymptotic potential, corresponding to the apparent water stream potential
in the moving reference frame of the body. As the name suggests,
the perturbation potential accounts for the perturbation effects that the hull presence has
on the flow field. Instead, the asymptotic potential describes the flow that would be observed
in the hull absence. As we will see in more detail in the numerical result sections, the
asymptotic potential is a known scalar function satisfying the Laplace equation and can
represent a steady stream velocity $\Ub_\infty = \nablab\phi_\infty$, a wave flow field such as the one
described by Airy potential, or a combination of both. 

The equations of motion that describe the velocity and pressure
fields $\ub(\xb,t)$ and $p(\xb,t)$  in the fluid region surrounding the
moving body are the incompressible Navier--Stokes equations. For a detailed derivation potential
flow equations applied to nonlinear water wave, and the corresponding boundary value problem, we
refer the interested readers to \cite{waveBem}. In the present discussion, we will report the main
results with focus on aspects that are relevant to obtaining a boundary value problem that, when
discretized using a boundary element method provides a unified framework for solving both
steady and unsteady ship wave problems. Under the aforementioned
assumptions, the continuity and momentum equation can be recast into the Laplace Equation for
the velocity potential and the unsteady Bernoulli equation, respectively. Under the assumption $\Delta\phi_\infty=0$
we have
\begin{align}
\label{eq:pot}\Delta \phi & =  0 & \text{in}\  \Omega(t),\ \forall t\geq 0\\
\label{eq:bern}\dsfrac{\Pa\phi}{\Pa t} + \dsfrac{\Pa\phi_\infty}{\Pa t} + \dsfrac{1}{2}|\nablab\phi_\infty+\nablab\phi|^2+\gb\cdot\xb+\dsfrac{p}{\rho} 
& =  C(t) &  \text{in}\  \Omega(t),\ \forall t\geq 0.
\end{align}
Here, $\rho$ indicates the --- constant --- density of the fluid and the reference frame gravity acceleration
vector $\gb$ is used in the corresponding gravity forces potential term. In particular, in the case of inertial reference
frame, $\gb=g\eb_z$, where $g=9.81\,\text{m}/\text{s}^2$ is the earth gravity acceleration. 
Since pressure only appears in Equation \eqref{eq:bern}, a typical approach in potential flow theory
is that of solving Equation \eqref{eq:pot} to obtain the perturbation potential, which is then
introduced in Bernoulli's equation to evaluate the pressure field. Thus, our governing equation is the Laplace
equation for the perturbation potential field, from which both velocity and pressure fields can be
recovered by means of Equations \eqref{eq:pert_pot} and \eqref{eq:bern}, respectively.
To obtain a well posed problem for the perturbation potential field, the Laplace equation
must be complemented by a suitable set of conditions on the domain
boundary $\partial\Omega(t)=\Gamma(t)=\Gamma^b\cup\Gamma^h\cup\Gamma^f(t)\cup\Gamma^i\cup\Gamma^\infty,\ \forall t\geq 0$. In this work,
we consider the free surface $\Gamma^f(t)$ as the only part of the boundary not fixed, and we will drop the explicit dependence
on time for the sake of simplicity.  We however point out that in different investigations, the model presented
would be able to consider also the motion of other boundaries --- such as for instance the bottom or the hull one --- that
are not considered here. On the
bottom boundary $\Gamma^b$ of the basin  (assumed horizontal and located at $z=-H$) we set a non penetration boundary condition, namely

\begin{equation}
\ub(\xb,t)\cdot\nb = (\nablab\phi_\infty(t)+\nablab\phi(\xb,t))\cdot\nb = 0  \qquad \text{on}\ \Gamma^b.
\end{equation}
We assume that $\nablab\phi_\infty(t,\xb)|_{z=-H}$ is tangent to $\Gamma^b$, which is the case for a uniform flow field or Airy potential. 
This results in the following
homogeneous Neumann boundary condition for the perturbation potential

\begin{equation}
\dsfrac{\Pa\phi(\xb,t)}{\Pa n} = \nablab\phi(\xb,t)\cdot\nb = 0  \qquad \text{on}\ \Gamma^b.
\end{equation}
A non penetration boundary condition is also used on the body surface $\Gamma^h$,namely

\begin{equation}
\ub(\xb,t)\cdot\nb = (\nablab\phi_\infty(t)+\nablab\phi(\xb,t))\cdot\nb = 0  \qquad \text{on}\ \Gamma^h.
\end{equation}

In this case,
no assumptions on the relative orientation of $\nablab\phi_\infty(t)$ and $\nb$ can be made. Thus, the
non homogeneous Neumann body boundary condition for the perturbation reads

\begin{equation}
\dsfrac{\Pa\phi(\xb,t)}{\Pa n} = \nablab\phi(\xb,t)\cdot\nb =  -\nablab\phi_\infty(t)\cdot\nb  \qquad \text{on}\ \Gamma^h.
\end{equation}
In this work, we will make use of a homogeneous Dirichlet boundary condition on the inflow boundary $\Gamma^i$ of the domain,
 translating the fact that the perturbation potential must fade at infinite distance from the body. Considering
Equation \eqref{eq:pert_pot} the Dirichlet condition, which reads 

\begin{equation}\label{eq::BC_inflow}
\phi(\xb,t) = 0  \qquad \text{on}\ \Gamma^i,
\end{equation}
directly imposes that on the inflow boundary the flow field must be coincident with the one characterized by the asymptotic potential.  In other
words, we are assuming that the perturbation of the asymptotic flow associated with the presence of the body is null
at infinite distance upstream from the body. This is an acceptable assumption, as in the three dimensional setup here considered
the waves and the velocity perturbations generated by a finite dimensional body get dispersed over a wider and wider area
as they travel at greater distance from their source and ultimately fade to zero value.

A further boundary condition is to be applied on the far field truncation boundary $\Gamma^\infty$ of the numerical domain.
Ideally, also such condition should be able to translate the fact that the perturbation potential must fade to zero at
great distance from the body causing the perturbation. At the same time, the boundary condition should be
neutral to water waves reaching the boundary, avoiding their reflection into the basin. Given the dispersive
nature of water gravity waves reaching the outer boundary of the domain, devising a wave absorbing boundary condition
working effectively across a wide range of wavelengths is a rather difficult task, especially in three dimensions.
For such a reason, there has been wide debate over the most effective form of the boundary condition to be applied at
the far field truncation boundary of the numerical domain, and many investigations have been carried
out on the subject (see for instance \cite{num_beach}). Discussing the strengths and weaknesses of each approach presented in the
literature is clearly beyond the scope of this work, in which we have made use of a simple homogeneous Neumann boundary
condition, namely

\begin{equation}
\dsfrac{\Pa\phi(\xb,t)}{\Pa n} = \nablab\phi(\xb,t)\cdot\nb = 0  \qquad \text{on}\ \Gamma^\infty.
\end{equation}
As such condition does not address the problem of waves reflection, it has been complemented by
the presence of a numerical damping zone --- or \emph{numerical beach} --- located immediately before the downstream
and upstream boundaries of the domain. While the general idea for such damping zone was originally presented
in \cite{num_beach}, it had to be adapted to the numerical setup used in the present work. Because the numerical beach
is implemented including an additional term in the free surface boundary condition, we will detail its
implementation in the next section, devoted to such boundary condition .  

The correct conditions to be applied on the free surface boundary $\Gamma^w$ to allow for an accurate
tracking of the water waves are in fact one of the most interesting and delicate modeling aspects of the application of 
potential flow theory to this kind of flows. In the next section, we will detail several possible choices, and discuss
their implications on the well-posedness of the mathematical problem on both steady and unsteady assumptions.

\subsection{Free surface boundary conditions}

The first thing to be pointed out about the free surface boundary $\Gamma^w$ is that its position is unknown
\emph{a priori}. For such a reason,  given the additional unknown of the problem, an additional boundary condition
must be specified on the free surface boundary. More specifically, the discussion on the correct treatment
of $\Gamma^w$ must not only address the problem of finding the most suitable conditions for the perturbation
potential $\phi$, but also the correct way to update the free surface position during the time integration of the problem. 

In this regard, we first focus  our attention on the governing equation for
the free surface position. We start assuming that the water free surface
elevation field $\eta$ is a single valued function of the horizontal coordinates $x$ and $y$, namely

\begin{equation} \label{eq:eta_one_valued}
z = \eta(x,y,t)  \qquad \text{on}\ \Gamma^f.
\end{equation}

It is quite clear that such an assumption limits the applicability of this model to cases in which no wave breaking
occurs. However, wave overturning might not only result in free surface making contact with itself, but
also in the presence of vortical flow regions.  Thus, we remark that considering the presence of breaking waves would require 
abandoning the potential flow model altogether. The development of a multi-model solver in which the potential flow equations
are interfaced with viscous models in the regions in which the flow is vortical is definitely an extremely interesting
research area. Yet, it again falls far beyond the scope of the present work, which instead aims at obtaining a
free surface potential flow model which is efficient and robust enough to be possibly interfaced with other models.

\subsubsection{Derivation of Lagrangian free surface boundary conditions}

So, the knowledge of the free surface elevation field $\eta(x,y,t)$ results in the complete description of the domain
shape $\Omega(t)$. $\eta(x,y,t)$ will then become one of the unknowns of the mathematical problem at hand. To write an
evolution equation for such a new variable, we move from assumption \eqref{eq:eta_one_valued} to obtain
a constraint $\mathcal{G}(x,y,z,t)$ which reads

\begin{equation} \label{eq:constraint}
\Gamma^f(t) = \{x \in \partial \Omega(t) : \mathcal{G}(x,t) =z-\eta(x,y,t) =0\}\ \forall t\geq 0.
\end{equation}
Taking the Lagrangian derivative (considering $f(\xb, t)$ we define $\dsfrac{D f}{D t}=\dsfrac{\Pa f}{\Pa t}+\dsfrac{\Pa \xb}{\Pa t}\cdot \nablab f = \dsfrac{\Pa f}{\Pa t}+\ub \cdot \nablab f $)
of $\mathcal{G}(x,y,z,t)$ we obtain

\begin{equation} \label{eq:constraint_lag_der}
\frac{D\mathcal{G}}{D t} = \frac{D z}{Dt} - \frac{D\eta}{Dt} = 0 \qquad \text{on}\ \Gamma^f,
\end{equation}
which results in 

\begin{equation} \label{eq:kin_fs_condition}
\frac{D\eta}{Dt} = \frac{\Pa \phi_\infty}{\Pa z}+\frac{\Pa \phi}{\Pa z} \qquad \text{on}\ \Gamma^f.
\end{equation}

Equation \eqref{eq:kin_fs_condition} is referred to as the Lagrangian form of the fully nonlinear
free surface \emph{kinematic} boundary condition. A similar boundary condition, used to
update the free surface values of the perturbation potential $\phi$, is obtained from
the manipulation of Bernoulli equation. Adding $\frac{1}{2}|\nablab\phi_\infty+\nablab\phi|^2$ to both sides
of Equation~\eqref{eq:bern} yields

\begin{equation} \label{eq:dyn_fs_condition_step_1}
\dsfrac{\Pa\phi}{\Pa t} + \dsfrac{\Pa\phi_\infty}{\Pa t}+ |\nablab\phi_\infty+\nablab\phi|^2+\gb\cdot\xb+\dsfrac{p_a}{\rho} 
=  C(t) + \dsfrac{1}{2}|\nablab\phi_\infty+\nablab\phi|^2\quad  \text{on}\  \Gamma^f,
\end{equation}
which, working out the computations and  making use of $\gb\cdot\xb=g\eta$ results in

\begin{equation} \label{eq:dyn_fs_condition_step_2}
\dsfrac{D\phi}{ Dt}  
=  C(t)  - \dsfrac{\Pa\phi_\infty}{\Pa t} - g\eta-\dsfrac{p_a}{\rho} + \dsfrac{1}{2}|\nablab\phi|^2 - \dsfrac{1}{2}|\nablab\phi_\infty|^2 \quad  \text{on}\  \Gamma^f.
\end{equation}

We now assume that the atmospheric pressure $p_a$ exerted by air on the water free surface is a constant
and uniform field, and that its value --- which is defined up to a constant --- is set to zero. The assumption
that the air pressure on water is a uniform field is a 
rather reasonable one, especially in presence of streamlined displacement vessels traveling at moderate speeds. In addition,
since the perturbation potential is assumed to fade for $|\xb|\rightarrow\infty$, it is possible to compute
that for a point on the free surface $C(t)=\dsfrac{\Pa\phi_\infty}{\Pa t}+\frac{1}{2}|\nablab\phi_\infty|^2+g\eta_\infty$. We must point out that
despite the fact that $\phi_\infty$ and $\eta_\infty$ depend on the position in which they are computed, because the far field potential
satisfies Bernoulli's equation the value of $C(t)$ is independent of the position. This is certainly true in the case in which the 
value of $\phi_\infty=\Ub_\infty(t)\cdot\xb$ represents a uniform stream velocity associated with flat free surface $\eta_\infty=0$.
This is also true for potentials representing an incident wave field, such as that of a monochromatic Airy wave
considered in the result cases. However, because the wave elevation field $\eta_\infty$ in the latter case is obtained making use of a
linearized free surface boundary condition, assuming that $C(t)$ is independent on the position is a good approximation only for
waves of small amplitude. Given all these considerations, the value of the final form of the fully nonlinear \emph{dynamic} free
surface boundary condition for the perturbation potential reads

\begin{equation} \label{eq:dyn_fs_condition}
\dsfrac{D\phi}{ Dt}  
=   \dsfrac{1}{2}|\nablab\phi|^2 - g(\eta-\eta_\infty)\quad  \text{on}\  \Gamma^f.
\end{equation}

The combined enforcement of Equations \eqref{eq:kin_fs_condition} and \eqref{eq:dyn_fs_condition} on the free
surface portion of the boundary $\Gamma^f$, effectively closes the mathematical problem, and allows for the
solution of a well posed boundary value problem \cite{grilli2001}.

\subsubsection{Eulerian time derivatives and steady state}\label{sec::eulerian_and_steady}
{

The main problem associated with the numerical resolution of the model represented by
Equations \eqref{eq:kin_fs_condition} and \eqref{eq:dyn_fs_condition} is related to the presence of Lagrangian time derivatives.
In fact, as already pointed out in presence of a main stream velocity past the hull such derivatives cause the
downstream drift of the mesh nodes upon resolution of the numerical problem. It is of course possible to make use of
the Lagrangian derivative definition to only include Eulerian
derivatives in the free surface boundary conditions (Equations \eqref{eq:kin_fs_condition} and \eqref{eq:dyn_fs_condition}), as follows

\begin{eqnarray}
\label{eq:cin_fs_condition_eul}
\dsfrac{\Pa\eta}{\Pa t} &=& \frac{\Pa \phi_\infty}{\Pa z}+\frac{\Pa \phi}{\Pa z} - (\nablab\phi+\nablab\phi_\infty)\cdot\nablab\eta\qquad \text{on}\ \Gamma^f \\   
 \label{eq:dyn_fs_condition_eul}
\dsfrac{\Pa\phi}{\Pa t}  
&=&   \dsfrac{1}{2}|\nablab\phi|^2 -(\nablab\phi+\nablab\phi_\infty)\cdot\nablab\phi- g(\eta-\eta_\infty)\quad  \text{on}\  \Gamma^f
\end{eqnarray}

In principle, as the latter equation only requires the evaluation of Eulerian time derivatives it should not pose
problems related to downstream drift of the computational mesh nodes. Unfortunately, the numerical evaluation of the Eulerian derivatives
on a moving boundary such as the free surface would pose several problems. 
For instance in case of a collocation scheme --- the most common choice in BEMs --- the collocation points on which the
solution vector is computed would not sit on the same spatial position for all the time steps needed for the time derivative
evaluation. 

Nonetheless, several researchers took advantage of the fact that if and when a steady solution is reached, the Eulerian
derivatives are null and can be eliminated from Equations \eqref{eq:cin_fs_condition_eul} and \eqref{eq:dyn_fs_condition_eul}.
From that starting point, a series of different methods have been developed for the solution of steady ship hydrodynamics
problems. As already pointed out in Section~\ref{sec:intro}, Raven \cite{ravenPhD1998} combined a non dimensional form of
Equations \ref{eq:cin_fs_condition_eul} and \ref{eq:dyn_fs_condition_eul} to obtain a single free surface boundary condition. In
the framework of an iterative method in which the free surface elevation $\eta$ is known from the previous iteration, such
combined free surface boundary condition is enforced in the Laplace problem for the perturbation potential. The potential solution
obtained is then used to update the free surface elevation and move to the ensuing iteration, until convergence is reached.
As briefly reported in Section~\ref{sec:intro}, Scullen~\cite{Scullen} employed a different and interesting approach in which a non
penetration Neumann boundary condition was
enforced on the free surface in the Laplace problem for the perturbation potential. The resulting potential solution was then
used to compute the new free surface elevation making use of a dynamic condition.
Scullen tested four different possible free surface dynamic conditions similar to
Equation \eqref{eq:dyn_fs_condition_eul} (imposing null pressure, null pressure Lagrangian derivative, and two different
combinations of the latter quantities), which all led to convergence to steady solutions. These efforts led to a series of
solvers that are able to obtain remarkably
good estimates of the wave resistance of a ship advancing at constant speed in calm water, and have been used with success
in the last three decades in the design process of a vast number of ships. But despite their success in solving problems
characterized by steady flows, these models cannot be easily adapted to solve time dependent problems. In fact, simply adding
back the Eulerian time derivatives disregarded in the steady case would not help, given the aforementioned problems in their computation
on moving domains.



\subsubsection{Semi-Lagrangian (or Arbitrary Lagrangian--Eulerian) free surface boundary conditions}\label{sec::semi-lag_BC}

To sum up, using Lagrangian derivatives in the  the fully nonlinear free surface boundary conditions results in computational
grid nodes which, even in presence of steady flow, never really settle for an equilibrium position but instead drift around
and eventually away from the surging body. Considering instead the alternative of resorting to Eulerian derivatives
in the free surface conditions, also results in an unsatisfactory result. In fact, Eulerian derivatives are not suited with
the presence of a domain with moving boundaries and moving meshes. As suggested
by several references (\cite{doneaetAlAle2004,Formaggia2009CardiovascularM,QUARTERONI20043}) on fluid dynamic
applications of the Finite Element Method (FEM) in deforming domains, the problem related to time derivatives
in the fully nonlinear boundary conditions
can be solved resorting to the Arbitrary Lagriangian Eulerian (ALE) formulation. Introducing an arbitrary grid velocity
field $\vb(\xb,t):\ \mathbb{R}^3\rightarrow\mathbb{R}^3$ allows for the definition of the following total derivative.
Given a generic scalar field $f(\xb,t)$ the ALE time derivative reads

\begin{equation} \label{eq:total_der}
\frac{\delta f}{\delta t} = \frac{\Pa f}{\Pa t} + \vb\cdot\nablab f \ .
\end{equation}
We remark that derivative $\frac{\delta}{\delta t}$ represents the time derivative of the desired scalar field,
computed following points moving according to the prescribed --- grid --- velocity field $\vb$.
{

Adding the term
$\vb\cdot\nablab \eta$ on both sides of Equation \eqref{eq:kin_fs_condition} results in

\begin{equation} \label{eq:kin_fs_condition_trans}
\frac{D\eta}{Dt} + \vb\cdot\nablab \eta = \frac{\Pa\eta}{\Pa t} + (\nablab\phi_\infty+\nablab\phi)\cdot\nablab \eta + \vb\cdot\nablab \eta = \frac{\Pa \phi_\infty}{\Pa z}+\frac{\Pa \phi}{\Pa z} + \vb\cdot\nablab \eta\qquad \text{on}\ \Gamma^f.
\end{equation}
Rearranging the terms leads to

\begin{equation} \label{eq:kin_fs_condition_semil}
\frac{\delta\eta}{\delta t} = \frac{\Pa \phi_\infty}{\Pa z}+\frac{\Pa \phi}{\Pa z} + (\vb-\nablab\phi_\infty-\nablab\phi)\cdot\nablab  \eta\qquad \text{on}\ \Gamma^f.
\end{equation}
With a similar treatment, Equation \eqref{eq:dyn_fs_condition} becomes

\begin{equation} \label{eq:dyn_fs_condition_semil}
\dsfrac{\delta\phi}{ \delta t}  
=  \dsfrac{1}{2}|\nablab\phi|^2  -  g(\eta-\eta_\infty)\quad + (\vb-\nablab\phi_\infty-\nablab\phi)\cdot\nablab \phi \ \ \ \text{on}\  \Gamma^f.
\end{equation}

Equations \eqref{eq:kin_fs_condition_semil} and \eqref{eq:dyn_fs_condition_semil} are the ALE formulation
of the Kinematic and Dynamic fully nonlinear free surface boundary conditions. In the naval engineering
literature they are typically referred to as fully nonlinear free surface boundary conditions written in
semi-Lagrangian form. They have been first introduced by Beck in \cite{beck1994}, which pointed out that
velocity field $\vb$ can be selected so as to avoid the downstream motion of the nodes. In particular
choosing

\begin{equation} \label{eq:v_grid}
\vb=v_x\eb_x+v_y\eb_y+\dsfrac{\delta\eta}{\delta t}\eb_z  \quad \text{on}\ \ \Gamma^f(t)
\end{equation}
allows for the computation of the vertical velocity $\frac{\delta\eta}{ \delta t}$ required to keep on the free
surface a point moving with arbitrary horizontal velocities $v_x$ and $v_y$.  As for the value of $\vb$ on all
the other boundaries, the imposed value is null, as in this work we are not considering possible motion of
hull, bottom or far field surfaces.

The grid velocity field $\vb$ is in principle defined on the whole three dimensional domain $\Omega$. However,
we must remark that in this work we a employ a Boundary Element Method for the spatial discretization of the boundary value
problem equations. For such a reason, only the values of $\vb$ on the free surface are needed. Should other 
codimension zero discretization methods such as the Finite Element method be used, it would be possible to resort to
an harmonic extension to obtain the values of $\vb$ in $\Omega$ based on the ones computed on its boundary. 
The latter procedure is commonly used in Fluid-Structure Interaction simulations to extend the ALE map
from the fluid boundary to the internal portion of the computational domain. For a more detailed description
of such a procedure and of the ALE approach in general, we once again refer the interested
reader to \cite{Formaggia2009CardiovascularM} (Chapter 3, Section 3.5 in particular).

At the numerical level, the initial boundary value problem described by these equations is solved by means
of an approach named Mixed Eulerian Lagrangian scheme. In such framework, a Dirichlet boundary condition
is prescribed on the free surface in the Laplace problem for the perturbation potential. Once the potential is
known from the Laplace problem solution, the resulting fluid velocity is first used in Equation~\eqref{eq:kin_fs_condition_semil}
to compute the new position of the free surface collocation points, and then used in \eqref{eq:dyn_fs_condition_semil}
to compute the potential on the displaced free surface. With a new free surface position and potential, a new
time step can be started to continue the computation.

The potential flow model
employing the ALE formulation of the fully nonlinear free surface boundary conditions reported in
Equations \eqref{eq:kin_fs_condition_semil} and \eqref{eq:dyn_fs_condition_semil} has been successfully
implemented in several contributions (\cite{kjellbergEtAl2011,scorpioPhD1997, waveBem}). Its main advantage
with respect to the model characterized by the Lagrangian boundary conditions, is that the arbitrary grid
velocity components $v_x$ and $v_y$ can be selected so as to avoid downstream drift of the grid nodes and retain mesh quality.
On the other hand,  a comparison of the ALE boundary conditions in Equations \eqref{eq:kin_fs_condition_semil}
and \eqref{eq:dyn_fs_condition_semil} with the Lagrangian counterparts in Equation~\eqref{eq:kin_fs_condition}
and \eqref{eq:dyn_fs_condition}, suggests that transport terms appear in the ALE formulation. Such transport terms,
in which the advection vector is the difference between the fluid velocity and the grid velocity, were in fact not present
in the Lagrangian boundary conditions. It is at this stage important to point out that whenever
the grid velocity $\vb$ becomes significantly different with respect to the fluid velocity $\nablab\phi_\infty+\nablab\phi$,
such a transport term can become dominant, leading to possible stability issues in the problem discretization.
Since keeping $v_x$ and $v_y$ values close to 0, even in presence of high stream velocity, is necessary
to avoid downstream drift of the grid nodes, a specific stabilization method will be needed to fully
exploit the advantages of the ALE formulation, as discussed in Section~\ref{sec::SUPG}.

The free surface model described by
Equations \eqref{eq:kin_fs_condition_semil} and \eqref{eq:dyn_fs_condition_semil}
was implemented in unsteady flow solvers featuring an implicit and stabilized version of the MEL time advancing scheme,
that proved capable of simulating the flow past ship hulls
without downstream drift of the free surface grid and without remeshing was required \cite{waveBem}. However,
despite such flow admits steady solutions, the regime solution could only be obtained using the unsteady formulation and
integrating the equations for a sufficiently large number of time steps. In fact, at the numerical level
the steady nonlinear problem obtained removing time derivatives from Equations \eqref{eq:kin_fs_condition_semil}
and \eqref{eq:dyn_fs_condition_semil} was never able to converge to meaningful solutions. In the authors estimation,
such failure to reach convergence could be caused by the choice of using a Dirichlet boundary condition on the free surface
in the boundary value problem for the perturbation potential. In several successfully converging steady flow
solvers in fact, the nonlinear iterations are started from a so called \emph{double body potential} solution
obtained imposing a non penetration Neumann boundary condition on both the hull surface and on the free surface
sitting in its flat undisturbed position. Such solution is relatively close to the final solution, as the fluid
dynamic field and is already accounting for both the presence of the hull and that of the free surface boundary,
and only accurate adjustment of the free surface elevation is required to finally reach convergence. That is why
steady solvers such as the one implemented by Scullen \cite{Scullen} impose a non penetration Neumann boundary condition
on the free surface in the perturbation potential Laplace problem, and then use a \emph{dynamic} condition
such as Equation \eqref{eq:dyn_fs_condition_eul} to compute the equilibrium position of the free surface.

As opposed to this, in the MEL framework a Dirichlet boundary condition is first applied on the free surface in the
Laplace problem for the perturbation potential. After this, a \emph{kinematic} condition such as
Equation \eqref{eq:kin_fs_condition_semil} condition is used to find the equilibrium position of the free surface
and take care of its displacement. Arguably, using a --- non penetration --- condition inspired on kinematic considerations
as free surface boundary condition of the Laplace problem for the perturbation potential, and a dynamic condition
to describe the shape and motion of the free surface appears to have more physical meaning. In fact, when
a MEL simulation is started from a double body potential flow, a kinematic condition such as Equation \eqref{eq:kin_fs_condition_semil}
is unable to lead to any change in the free surface elevation. So, instead of being started from such convenient
initial guess, the steady solver obtained dropping the time derivatives in Equations \eqref{eq:kin_fs_condition_semil}
and \eqref{eq:dyn_fs_condition_semil} can only start from a null solution, which is not as close to the target one.
For such a reason, in this work we replaced the kinematic condition expressed by Equation \eqref{eq:kin_fs_condition_semil},
with a non penetration Neumann boundary condition which is used on the free surface in the Laplace problem
for the perturbation potential. As a consequence, the free surface motion will be determined based on the
dynamic condition given by Equation \eqref{eq:dyn_fs_condition_semil}. Along with writing the corresponding
boundary value problem, in the next section we will prove that under the assumption that $\eta$ is a 
Cartesian function of $x$ and $y$, the non penetration boundary condition used in the present work can indeed be derived
by the kinematic condition used in previous works. We must remark that the latter prove is not strictly necessary, since
there is no requirement that different physical models for the same phenomenon should result in the same set of
equations at the mathematical level. However, in the present case the fact that the Neumann boundary condition can be derived
from the nonlinear kinematic free surface boundary condition might indeed explain why the results obtained with the two models
are in close agreement.

\subsubsection{ Alternative kinematic free surface boundary condition formulation}

It is quite interesting to point out that, by a physical standpoint, the fully nonlinear kinematic
free surface boundary condition in Equation \eqref{eq:kin_fs_condition} represents the requirement
that a material point on the free surface will remain on the free surface, which indeed is a stream
surface for the fluid velocity field. In fact, the equations states that the Lagrangian time
derivative of the coordinates of any point on the free surface, is equal to the local Eulerian velocity.
Clearly, this consideration applies to the semi-Lagrangian boundary condition in Equation \eqref{eq:kin_fs_condition} too,
as it is derived by it Lagrangian counterpart. So, since the free surface is a stream surface for the
velocity field, intuition suggests that a simple non penetration boundary condition can be applied
on such a boundary portion, rather than Equations \eqref{eq:kin_fs_condition} and \eqref{eq:kin_fs_condition_semil}.
To prove that a non penetration condition can be derived from \eqref{eq:kin_fs_condition_semil} we start
considering the free surface tangent plane equation. At a generic point $\xb\in \Gamma^f$, in which
the outward unit normal vector to the free surface boundary $\Gamma^f$ is $\nb = n_x\eb_x+n_y\eb_y+n_z\eb_z$,
the implicit equation of the tangent plane reads

\begin{equation} \label{eq:impl_tang_plane}
n_x x+n_y y+n_z z=c
\end{equation}
from which we can obtain the Cartesian equation of the tangent plane along direction $z$

\begin{equation} \label{eq:cart_tang_plane}
 z=-\frac{n_x}{n_z} x-\frac{n_y}{n_z} y+\frac{c}{n_z}.
\end{equation}
By its very definition, the  slope of the tangent plane at point $\xb$ coincides with
that of the free surface at the same point, which gives us a simple way to compute the free surface
elevation gradient as
\begin{equation} \label{eq:eta_grad}
\nablab\eta=-\frac{n_x}{n_z}\eb_x-\frac{n_y}{n_z}\eb_y.
\end{equation}
Substituting the latter gradient representation in Equation \eqref{eq:kin_fs_condition_semil}, we have

\begin{eqnarray} \label{eq:kin_fs_condition_semil_modif}
n_z \frac{\delta\eta}{\delta t} &=& n_z\nablab\phi_\infty\cdot\eb_z+n_z\frac{\Pa \phi}{\Pa z} - n_x\vb_x -  n_y\vb_y +\nonumber \\
&& n_x\nablab\phi_\infty\cdot\eb_x+ n_y\nablab\phi_\infty\cdot\eb_y+n_x\frac{\Pa\phi}{\Pa x}+n_y\frac{\Pa\phi}{\Pa y}  \qquad \text{on}\ \Gamma^f.
\end{eqnarray}
Rearranging terms and making use of  the definition of $\vb$ in Equation~\eqref{eq:v_grid} we obtain 

\begin{equation} \label{eq:kin_non_pen_ALE}
\vb\cdot\nb= \left(\nablab\phi+\nablab\phi_\infty\right)\cdot \nb  \qquad \text{on}\ \Gamma^f,
\end{equation}
which is the non penetration boundary condition for the perturbation potential in
presence of a moving boundary. In fact, it states that the normal component of the fluid velocity
must be equal to the normal component of the local boundary velocity. For our purposes,
Equation \eqref{eq:kin_non_pen_ALE} can be finally recast in the form of the following non homogeneous
Neumann boundary condition

\begin{equation} \label{eq:kin_fs_condition_non_pen}
\frac{\Pa \phi}{\Pa n} = \nablab\phi\cdot\nb= \left(\vb-\nablab\phi_\infty\right)\cdot \nb  \qquad \text{on}\ \Gamma^f.
\end{equation}

A further look at Equation~\eqref{eq:kin_fs_condition_non_pen} suggests that it represents an ALE generalization
of a Neumann non penetration boundary condition for the case of moving boundaries. For this reason, we can in principle
apply a similar, more general, condition also to other boundaries as the hull and bottom surface, in which non
penetration is enforced. That is why we will now write the final form of the boundary value problem solved in this
work making use of this form of the Neumann condition on all the non penetration boundaries. A specific choice
of boundary velocity $\vb$ will then be made for each of these boundaries.

\subsection{The final boundary value problem}

So, after introducing the governing equation along with the boundary conditions, here is the complete
boundary value problem considered in the present work. We are looking for $\phi$ and $\eta$ such that

\begin{subequations}\label{sys:full_BVP}
\begin{empheq}[left = \empheqlbrace]{align}
\Delta \phi &= 0 & \text{in} &\ \Omega(t),\ \forall t\geq 0 \\
\label{eq:neumann_ALE}
\frac{\Pa \phi}{\Pa n} &= \left(\vb-\nablab\phi_\infty\right)\cdot \nb  & \text{on}&\ \Gamma^b \cup \Gamma^h \cup \Gamma^f(t),\ \forall t\geq 0 \\
\frac{\Pa \phi}{\Pa n} &= 0  & \text{on}&\ \Gamma^\infty,\ \forall t\geq 0 \\
\phi &= 0 &  \text{on}&\ \Gamma^i,\ \forall t\geq 0 \\
\label{eq:dynamic_ALE}
\nonumber\frac{\delta \phi}{\delta t} &= \dsfrac{1}{2}|\nablab\phi|^2 - g(\eta-\eta_\infty)+& \\
                             & \ \ \    (\vb-\nablab\phi_\infty-\nablab\phi)\cdot\nablab \phi  & \text{on}&\ \Gamma^f(t),\ \forall t\geq 0,
\end{empheq} 
\end{subequations}

with

\begin{eqnarray}
\vb=v_x\eb_x+v_y\eb_y+\dsfrac{\delta\eta}{\delta t}\eb_z  & \text{on}&\ \Gamma^f(t) \\
\vb = 0 &\text{on}& \ \Gamma^b \cup \Gamma^h.
\end{eqnarray}

Thus, in this work the grid velocity $\vb$ appearing in the ALE formulation of the Neumann condition imposed on
$\Gamma_b$ and $\Gamma_h$ is null. This means that the computational mesh nodes on both the hull and the bottom surfaces
are idle. On the other hand, the vertical position $\eta$ of the computational nodes on the free surface
boundary $\Gamma^f(t)$ is an unknown. On such boundary $\vb$ is computed from $\delta\eta/\delta t$ according to Equation~\eqref{eq:v_grid}.

We point out that using the ALE formulation of non penetration Neumann boundary conditions represented by Equation \eqref{eq:neumann_ALE}
allows in principle to simulate the unsteady flow driven by the non stationary motion of the hull surface or on the bottom boundary. In such
case, the ALE velocity $\vb(t)$ prescribed in correspondence with such boundaries would be different from zero. The
present formulation is also suited for free surface interaction problem, in which the boundary grid velocity field
on $\Gamma_b$ and/or $\Gamma_h$ is an additional unknown of the problem. In such case, additional equations for the
dynamics of the hull surface and/or bottom boundaries must be added to close the mathematical problem. For instance,
a set of three dimensional rigid body dynamics equations could be used to compute the motion of the hull under the
action of the hydrodynamic forces as in \cite{molaIsope2016}, and the resulting grid nodes velocities could be
introduced in Equation \eqref{eq:neumann_ALE} to provide a suitable interface with the potential fluid solver.

\subsection{The numerical beach damping term}

A draw back of using an homogeneous Neumann boundary conditions for the vertical far field boundary condition is that it reflects energy back in the computational domain. We use an absorbing beach technique, in which we add an artificial damping
region away from the hull, used to absorb the wave energy. A damping term can be seen as an additional pressure $P$ acting
on the free surface. The resulting modified dynamic free surface boundary condition reads

\begin{equation} \label{eq:dyn_fs_condition_semil_beach}
\dsfrac{\delta\phi}{ \delta t}  
=  \dsfrac{1}{2}|\nablab\phi|^2  -  g(\eta-\eta_\infty)\quad + (\vb-\nablab\phi_\infty-\nablab\phi)\cdot\nablab \phi -\frac{P}{\rho} \ \ \ \text{on}\  \Gamma^f.
\end{equation}

The damping pressure $P$ used in this work is constructed as

\begin{equation} \label{eq:damping}
P=\mu\left(\frac{\delta\eta}{\delta t}-\frac{\delta\eta_\infty}{\delta t}\right).
\end{equation}

Ideally, this choice is able to dissipate any vertical velocity of the grid nodes not equal
to $\eta_\infty$, i.e. the one associated with the far field potential $\phi_\infty$.
We remark that in case the far field potential is associated with a steady stream and a null free surface elevation,
this damping term dissipates any vertical velocity of the free surface nodes.
Numerical observations suggest instead that the steady simulations do not require damping in
order to remain stable. That is why the damping strategy proposed is only active in the unsteady cases.
Clearly, the damping pressure is acting only in proximity of the inflow and outflow boundaries
of the domain. Thus, coefficient $\mu$ reads

\begin{equation} \label{eq:coeff}
\mu=\left(\frac{\max(|x|-x_d)}{L_{d}}\right)^2
\end{equation}
in which $x_d>0,\ L_d>0$ are the distance from the origin at which the damping zone starts acting and
the overall length of the damping zone, respectively.

\section{Boundary value problem discretization}

The literature on potential flow solvers with fully nonlinear free surface boundary conditions,
suggests that the most common way to integrate over time and space a time dependent boundary value problem such as
that in System \eqref{sys:full_BVP}, is the Mixed Eulerian-Lagrangian (MEL) scheme \cite{longuet-higginsCokeletMEL76}.
In such a framework,
at each time step a Laplace boundary value problem with Dirichlet boundary condition on the free surface is solved.
After this, the kinematic boundary condition is time integrated to obtain the new position of the free surface, and
the same is done with the  dynamic boundary condition so as to obtain the new value of the free surface potential.
With the new values of the free surface position and potential, a new time step can be processed.
As is made clear by its name, the original MEL procedure was designed to work with the Lagrangian form of free
surface boundary conditions. Yet a MEL-like algorithm can be also applied to semi-Lagrangian free surface conditions,
as seen in \cite{scorpioPhD1997}. Unfortunately, MEL cannot be applied to the free surface boundary condition formulation
adopted in this work. In fact, the kinematic semi-Lagrangian free surface boundary condition has been here replaced
with a non homogeneus Neumann boundary condition on such boundary, which makes MEL not applicable. Thus, we
have resorted to a different approach originally presented in \cite{waveBem}, which consists in carrying out the
spatial discretization of the governing problem,  to obtain a system of Differential Algebraic Equations (DAE).
In fact, the Boundary Element Method (BEM)  discretization of the Laplace
equation results in a set of algebraic equations, while the Finite Element Method (FEM) discretization of the
free surface boundary conditions leads to a set of differential equations. The DAE combining these different
equations is then solved by means of a Backward Difference Formula (BDF) time integration scheme. The next subsections
will present details of the BEM and FEM used to discretize over space the Laplace equation and the dynamic free surface
boundary condition, respectively. In addition, we will provide a description of the BDF time advancing scheme
used to solve the resulting DAE system.

\subsection{Laplace equation discretization based on Boundary Element Method}

In this work, we make use of the Boundary Element Method (BEM) discretization method for the
spatial discretization of the governing boundary value problem equations. In the context
of fully nonlinear free surface potential flow simulation, this is quite a common choice. We must
however remark that the Laplace equation for the velocity potential can be also discretized by means of
the Finite Element Method (FEM). In this regard, noteworthy works have been carried out by
Ma and Yan \cite{city4324} and more recently by Bermudez et Al. \cite{alfredo_bermudez_2021_4638166}.
At a first glance, it would appear that the most important advantage of BEM compared to FEM 
is the reduced number of unknowns associated with the
codimension one grid. Unfortunately, in the practice such an advantage is typically offset by the presence of 
a dense resolution matrix in the discretized algebraic system. Yet, there are other advantages of BEM
that made us favor it over FEM. In particular, in the context of the present physical problem, where
moving boundaries such as the free surface are present, the codimension one grids required by BEM are much 
easier to generate, deform and manage without significant quality drop.

As already pointed out,
two free surface boundary conditions appear in System~\eqref{sys:full_BVP}. In particular,
Equation~\eqref{eq:dynamic_ALE} is an additional boundary condition needed to determine the
value of the supplementary unknown $\eta$. In the present section, we will discuss the BEM discretization
of the Laplace System~\eqref{sys:full_BVP} devoid of Equation~\eqref{eq:dynamic_ALE}.
The specifics of the numerical discretization of the latter boundary condition will be reported in
Section~\ref{sec::SUPG}. The algorithm for the simultaneous solution of the two sub-problems will be then
described in Section~\ref{sec:neumann} and Section~\ref{sec:time_scheme}.
For the BEM discretization we use the same formalism presented in~\cite{Giuliani2015,brebbia,pi-BEM}, and here refer to a
bounded open domain $\Omega$ with Lipschitz boundary $\Gamma = \partial \Omega$. In such a framework, System \eqref{sys:full_BVP} -- devoid
of Equation~\eqref{eq:dynamic_ALE} --- is recast into the generic Laplace problem 

\begin{subequations}
\begin{alignat}{2}
  \label{laplaceint}
  -\Delta\phi &= 0  && \text{in }  \Omega\\
  \dfrac{\partial \phi}{\partial n} & = f_N(\mathbf{x}) \qquad && \text{on } \Gamma_{N} \\
   \phi&=f_D(\mathbf{x}) &&\text{on } \Gamma_{D},
 \end{alignat}
\end{subequations}
 where Dirichlet and Neumann boundary conditions are imposed on the portions $\Gamma_D$, and $\Gamma_N$ of $\partial \Omega$. We require that $\overline{\Gamma_D \bigcup \Gamma_N}=\partial \Omega$, $\Gamma_D\bigcap\Gamma_N=\emptyset$, we notice that $\Gamma_D \neq \emptyset$ ensures solution uniqueness.

\subsubsection{Boundary integral formulation}

To rewrite \eqref{laplaceint} as a Boundary Integral Equation (BIE) we make use of a \emph{fundamental solution} (or Green's function) of Laplace equation. More specifically, in this work we employ the free space Green's function

\begin{equation*}
G(\yb-\xb) = \frac{1}{4\pi|\yb-\xb|},
\end{equation*}
which is the distributional solution of 

\begin{eqnarray*}
-\Delta G(\yb-\xb)&=&\delta(\xb)  \quad \text{in}\ \mathbb{R}^3\\
\lim_{|\yb|\rightarrow \infty} G(\yb-\xb)& = &0,
\end{eqnarray*}
where $\yb\in\mathbb{R}^3$ is a generic point, and
 $\xb\in\mathbb{R}^3$ is the center of the Dirac delta distribution $\delta(\xb)$.
If we select $\xb$ to be inside $\Omega$, use the deﬁning property of
the Dirac delta and the second Green identity, we obtain

\begin{equation}
\label{integral}
  \phi(\mathbf{x}) =\int_{\Gamma}G(\mathbf{x}-\mathbf{y})\frac{\partial \phi}{\partial n}(\mathbf{ y})\d s_y
-\int_{\Gamma}\phi(\mathbf{ y}){\nablab G(\mathbf{x}-\mathbf{y})\cdot\nb(\yb)}\d s_y
\ \forall\mathbf{x}\in\Omega,
\end{equation}
{
where in this case $\yb\in\Gamma$ is a generic integration point on the domain boundary --- as indicated by the subscript
in differential $\d s_y$ --- and $\nb(\yb)$ is the outward unit normal vector to boundary $\Gamma$.
}

From~\eqref{integral} we notice that if the solution and its normal derivative on the boundary $\Gamma$ are known then the potential $\phi$ can be computed in any point of the domain. 
Considering the trace of~\eqref{integral} we can write the boundary integral form of the original problem as 
\begin{subequations} \label{BoundIntFormul}
\begin{align}
  \alpha(\mathbf{x})\phi(\mathbf{x})& =  \resizebox{.55\hsize}{!}{$\int_{\Gamma}G(\mathbf{x}-\mathbf{y})\frac{\partial \phi}{\partial n}(\mathbf{ y})\d s_y
-\int^{PV}_{\Gamma}\phi(\mathbf{ y}){\nablab G(\mathbf{x}-\mathbf{y})\cdot\nb(\yb)}\d s_y$}
&&\text{on }\Gamma \label{Laplacesetnbc}\\
  \dfrac{\partial \phi}{\partial n} & = f_N(\mathbf{ y}) && \text{on } \Gamma_{N}\label{Neumannsetnbc} \\
   \phi&=f_D(\mathbf{ y}) &&\text{on } \Gamma_{D}\label{Dirichletsetnbc},
\end{align}
\end{subequations}
where we have considered the Cauchy Principal Value (CPV) of the second integral. 
Then we write explicitly the boundary conditions~\eqref{laplaceint} in~\eqref{BoundIntFormul} obtaining
\begin{equation}
\begin{split} \label{BoundIntFormul-2}
 \chi_{\Gamma_N} (\mathbf{x}) \alpha(\mathbf{x})\phi(\mathbf{x}) -  &
 \int_{\Gamma_{D}}G(\mathbf{x}-\mathbf{y})\frac{\partial \phi}{\partial n}(\mathbf{ y})\d s_y\\
+&\int^{PV}_{\Gamma_N}\phi(\mathbf{ y}){\nablab G(\mathbf{x}-\mathbf{y})\cdot\nb(\yb)}\d s_y =\\
 -& \chi_{\Gamma_D} (\mathbf{x}) \alpha(\mathbf{x}) f_D(\mathbf{x}) + \int_{\Gamma_N}G(\mathbf{x}-\mathbf{y}) f_N(\mathbf{ y}) \d s_y \\
-&\int^{PV}_{\Gamma_D}f_D(\mathbf{ y}){\nablab G(\mathbf{x}-\mathbf{y})\cdot\nb(\yb)}\d s_y.
\end{split}
\end{equation}

We remark that $\alpha(\mathbf{x})$ is obtained from the CPV
evaluation of the integral involving the derivative of the { Green's function}, usually it represents the fraction of
solid angle of the domain $\Omega$ seen from the boundary
point $\mathbf{x}$. We use the generic characteristic function $\chi_A$ (which is one if $\mathbf{x}\in A$  and zero otherwise)
to split the term $\alpha(\mathbf{x})$ between Neumann and Dirichlet boundaries.

\subsubsection{Discretisation}
\label{sec:discretisation}
The numerical discretization of~\eqref{BoundIntFormul-2} leads to a real Boundary Element Method (BEM). The resolution of a BEM requires the discretization of the unknowns using functional spaces defined on a Lipschitz boundary. We address this problem introducing suitable discretizations for the { Neumann boundary unknown $\phi$ and for the Dirichlet boundary unknown $\phins$. Such discretizations are based
on standard Lagrangian finite element spaces defined on $\Gamma$}. We use the same functional space to describe the geometry, this setting is often referred to as Isoparametric BEM.

We define the computational mesh as a quadrilateral decomposition
${\Gamma}_h$ of the boundary $\Gamma$. We require that two cells $K,
K'$ of the mesh only intersects on common edges or vertices, and that there exist a mapping from a reference cell $\hat K$ to $K$ whose Jacobian is uniformly bounded away from zero for all cells $K$. 
{
To ease mesh generation, the simulation tool developed allows for the definition of a very coarse grid, which is then automatically refined on the user prescribed geometry up to the desired level of refinement. Following~\cite{pi-BEM-repo, dealII92} an interface to CAD files --- which are the most common tool to define arbitrary geometrical descriptions~\cite{heltaiEtAlACMTOMS2021,BangerthHeisterHeltai-2016-b,ArndtBangerthDavydov-2017-a} --- is used to specify the desired geometry of the hull. This feature has been employed in ship-wave simulations through BEM, \cite{waveBem,Mola2014,MolaHeltaiDeSimone-2017-a}, and~\cite{pi-BEM} presents an example of an aeronautics-like NACA wing
shape. As will be shown in Section~\ref{sec:test_case}, such a feature has been used in this work to refine an initial coarse computational
grid on the CAD surface of the spheroid geometry used for the numerical tests.
}

If $\phi$ and $\frac{\Pa \phi}{\Pa n}$ must lie in the spaces $V$ and $Q$, defined as
\begin{equation*}
\begin{aligned}
V &:= \left\{ \phi \in H^{\frac{1}{2}}(\Gamma)   \right\}   \\
Q &:= \left\{  \gamma \in H^{-\frac{1}{2}}(\Gamma)\right\},
\end{aligned}
\end{equation*}
where $\Gamma=\partial\Omega$, then the integrals in equation \eqref{Laplacesetnbc} are bounded.
$H^{\frac{1}{2}}(\Gamma)$ is the space of traces on $\Gamma$ of functions in $H^{1}(\Omega)$, while $H^{-\frac{1}{2}}(\Gamma)$ is its dual space. 
We construct the discretized spaces $V_h$ and $Q_h$  as conforming finite dimensional subspaces of $V$ and $Q$ respectively,
  \begin{subequations}\label{eq:shape_functs}
    \begin{align}
      V_h &:= \left\{ \phi_h \in L^2(\Gamma_h)  :  \phi_{h|K} \in \mathcal{Q}^{r}(K), \,K \in \Gamma_h   \right\}  \equiv \text{span}\{\psi_i\}_{i=1}^{N_V}  \\
      Q_h &:= \left\{ \gamma_h \in L^2(\Gamma_h)   :  \gamma_{h|K} \in \mathcal{Q}^{r}(K), \,K \in  \Gamma_h \right\} \equiv \text{span}\{\psi_i\}_{i=1}^{N_V} ,
    \end{align}
\end{subequations}
where  $\mathcal{Q}^{r}(K)$ is the space of polynomials of order $r$ in each coordinate direction.  In
principle these two spaces can be built independently, but in this work
we made use of the same Finite Element discretisation for both the primal
and the dual unknown, i.e., $V_h = Q_h \equiv \text{span}\{\psi_i\}_{i=1}^{N_V}$, $\psi_i$ being the shape function associated with the $i$-th degree of
freedom of the discretized space, for a more detailed analysis see \cite{GiulianiPoliMIThesis}. Following~\cite{pi-BEM} we use iso-parametric discretisations based on standard $Q_N$ Lagrangian finite elements, and
by collocating the support points of the geometry patches directly on the CAD surfaces. This work only reports results obtained with bi-linear
elements ($r=1$). The use of higher order bi-quadratic or bi-cubic elements has indeed been attempted, but it results in less robust
simulations, which typically fail to converge after a small number of adaptive refinement cycles are executed. Given the satisfactory convergence results obtained
in~\cite{pi-BEM} with higher order elements on several Laplace problems, such stability issues have been attributed to the discretization of the free
surface conditions that will be discussed in Section \ref{sec::SUPG}. As an alternative, it is possible that initial guess solutions that lead to convergence of
the nonlinear free surface boundary value problem with bi-linear elements, do not lead to convergence with higher order ones. Attempts to provide
higher order discretizations with solutions initial guesses obtained with bi-linear elements will be carried out in the near future.

The generic elements of the discretized spaces read

\begin{equation}
\phi_h(\xb)=\sum_{j=1}^{N_V}\hat{\phi}_j \psi_j(\xb), \qquad  \gamma_h(\xb)=\sum_{j=1}^{N_V}\hat{\gamma}_j \psi_j(\xb), 
\end{equation}
where $\hat{\phib},\hat{\gammab}$ represent the value at each collocation point of potential and potential normal derivative, respectively, and $N_V$ represents
the overall number of degrees of freedom of the discretized space. Finally,
we use the double nodes technique, \cite{grilli2001}, to ensure accuracy in the resolution of the BEM even when sharp edges are present. 

The collocation method is a common resolution technique for a BEM since it does not require any additional integration of~\eqref{BoundIntFormul-2}. For a deeper analysis of the accuracy of this setting the reader is referred to~\cite{pi-BEM}.
Collocating \eqref{Laplacesetnbc} produces the linear system 
\begin{equation}
(\alpha+N) \hat{\phib} - D \hat{\gammab} = 0,
\label{LaplaceLinearSystem}
\end{equation}
where
\begin{itemize}
\item{$\alpha$ is a diagonal matrix with the values
    $\alpha(x_i)$, where $x_i$ represents the i-th collocation point;}
\item{$N_{ij}=\sum_{k=1}^K \sum_q^{N_q} \dfrac{\partial G}{\partial n}(x_i
    - x_q) \psi_q^j J^k $, where $\hat{K}$ represents the reference
  cell and $J^k$ is the determinant of the first fundamental form for each panel k;}
\item{$D_{ij}=\sum_{k=1}^K \sum_q^{N_q} G(x_i
    - x_q) \psi_q^j J^k $.}
\end{itemize}

{

 When the collocation point lies inside the
cell where we are integrating we use bidimensional Lachat Watson quadrature
formulas to treat singular kernel integrals, see \cite{lachatWatson}. In any other case,
we make use of standard Gauss integration rules. 
}

\subsubsection{Numerical Implementation}
We use~\cite{pi-BEM-repo} as a backbone library for our work. In particular we use High Performance Computing libraries as \deal \cite{dealII92} and Trilinos~\cite{Heroux2005} to split the computational load between different processors and to tackle linear algebra. We achieve multicore parallelism using  Intel Threading Building Block (TBB)~\cite{Reinders2007}. A similar combination has been successfully applied to achieve high computational efficiency in fluid dynamics, as demonstrated in ASPECT~\cite{Kronbichler2012}. We remark that our BEM implementation greatly benefits from the distributed memory parallelism , due to the structure of the matrix assembling procedures see \cite{pi-BEM} for more details.

\subsection{Dynamic free surface boundary condition spatial discretization}\label{sec::SUPG}

Following the procedure outlined in \cite{waveBem}, to tackle the numerical discretization of Equation { \eqref{eq:neumann_ALE}}, we resort to its weak form, which reads

\begin{eqnarray}
\nonumber
\left(\dsfrac{\delta\phi}{ \delta t},\wb\right)  
&=&  \left(\dsfrac{1}{2}|\nablab\phi|^2  -  g(\eta-\eta_\infty)\quad + (\vb-\nablab\phi_\infty-\nablab\phi)\cdot\nablab {\phi} ,\wb\right) \\
\label{eq:dyn_fs_condition_semil_weak}
&=&\left(b^{\dot{\phi}},\wb\right).
\end{eqnarray}
Here, $\wb\in V$ is a test function and the notation

\begin{equation}
  \begin{split}
    \left(a,b\right)_w &= \int_{\Gamma^w(t)}ab\ \d\Gamma\\
    \left(a,b\right) &= \int_{\Gamma(t)}ab\ \d\Gamma
  \end{split}
\end{equation}
indicates a scalar product in the space $L^2(\Gamma)$. The discretization of Equation~\eqref{eq:dyn_fs_condition_semil_weak} is carried
out by means of a Galerkin Finite Element Method (FEM) based on the shape functions defined in Equation \eqref{eq:shape_functs}, and
results in the following system of algebraic equations

\begin{equation}\label{eq:L2_proj}
M\hat{\dot{\phib}} = \bb^{\dot{\phi}}
\end{equation}
where
\begin{itemize}
\item $M$ is a sparse mass matrix, the entries of which are given by
      $M_{ij}=\left(\psi_j,\psi_i\right)$;
\item $\bb^{\dot{\phi}}$ is a right hand side vector, with entries given by
      $\bb^{\dot{\phi}}_{i}=\left(b^{\dot{\phi}},\psi_i\right)$;
\item the entries of vector $\hat{\dot{\phi}}$ represent the nodal values
      of potential ALE time derivative $\frac{\delta \phi}{\delta t}$. 
\end{itemize}
There are several advantages associated with such an $L^2$ projection approach. First, it
avoids the evaluation of the potential gradients and surface normal vectors in correspondence
with the free surface collocation points, where such quantities are not single valued. In fact,
the integrals appearing in the weak formulation only require  the right hand side of
Equation~\eqref{eq:dyn_fs_condition_semil} to be evaluated on the numerical integration
scheme quadrature nodes, which fall within each quadrilateral cell. At such location, the
potential gradients and surface normal vectors are single valued, which results in an
accurate spatial integration scheme. An additional advantage is that matrix $M$ is
sparse, so assembling it only leads to a modest computational overhead with
respect to only assembling the BEM matrix. Finally, a further advantage of the $L^2$
projection approach, is that Equation~\eqref{eq:L2_proj} can be readily modified to include
stabilization terms able to avoid the dominant transport instabilities occurring for
high stream velocity, and discussed in Section~\ref{sec::semi-lag_BC}. As in \cite{waveBem},
we make use of a Streamwise Upwind Petrov--Galerkin (SUPG) stabilization
(for more detail, see \cite{hughes1979,tezduyar2003computation})
strategy to suppress
free surface instabilities that initially lead to in saw-tooth shaped free surface and
eventually result in simulation blow up. The SUPG stabilization consists in replacing
the plain $L^2$ projection in System~\eqref{eq:L2_proj} with the weighted projection

\begin{equation}\label{eq:L2_proj_SUPG}
\widetilde{M}\hat{\dot{\phib}} = \widetilde{\bb}^{\dot{\phi}}
\end{equation}
where
\begin{itemize}
\item The entries of $\widetilde{M}$ are given by
      $\widetilde{M}_{ij}=\left(\psi_j,\psi_i+\db\cdot\nablab_s\psi_i\right)$;
\item The entries of the right hand side vector $\widetilde{\bb}^{\dot{\phi}}$ are given by
      $\widetilde{\bb}^{\dot{\phi}}_{i}=\left(b^{\dot{\phi}},\psi_i+\db\cdot\nablab_s\psi_i\right)$;
\item Vector $\db=\tau\frac{\vb-\Ub_\infty-\nablab\phi}{|\vb-\Ub_\infty-\nablab\phi|}$ is aligned
      with the local velocity direction, with $\tau$ being a scalar coefficient proportional
      to the local mesh size. 
\end{itemize}

\subsection{Neumann boundary conditions spatial discretization}\label{sec:neumann}

Taking a look at boundary condition \eqref{eq:neumann_ALE}, we can immediately notice that both the normal vector $\nb$
and the values of the free surface nodes velocity $\vb$ appearing on the right hand side are not single valued. Thus,
we again resort to writing such equation in its weak form, which reads

\begin{equation}
\label{eq:neu_condition_weak_2}
\left(\dsfrac{\Pa\phi}{ \Pa n},\wb\right)  
=  \left( (\vb-\nablab\phi_\infty)\cdot\nb ,\wb\right) =\left(b^{\phi_n},\wb\right).
\end{equation}

The discretization of Equation~\eqref{eq:neu_condition_weak_2} results in the following
system of algebraic equations

\begin{equation}\label{eq:L2_proj_neu}
M\hat{\gammab} = \bb^{\phi_n}
\end{equation}
where
\begin{itemize}
\item $M$ is a sparse mass matrix, the entries of which are given by
      $M_{ij}=\left(\psi_j,\psi_i\right)$;
\item $\bb^{\phi_n}$ is a right hand side vector, with entries given by
      $\bb^{\phi_n}_{i}=\left(b^{\phi_n},\psi_i\right)$.
\end{itemize}

\subsection{Time advancing scheme}
\label{sec:time_scheme}

The spatially discretized resolution system can be recast in the following form

\begin{equation}\label{eq:DAE}
\Fb(\dot{\yb},\yb,t)=0
\end{equation}
where
\begin{equation}\label{eq:DAE_sol_vect}
\yb(t)=
\left\{
\begin{array}{c}
\{\hat{\phib}\} \\
\{\hat{\gammab}\} \\
\{\zb\} 
\end{array}
\right\}
=
\left\{
\begin{array}{c}
\left\{ \begin{array}{c}\hat{\phib}_D \\ \hat{\phib}_N \end{array}\right\} \\
\left\{ \begin{array}{c}\hat{\gammab}_D \\ \hat{\gammab}_N \end{array}\right\} \\
\left\{ \begin{array}{c}\zb_{FS} \\ \zb_{BA} \end{array}\right\} 
\end{array}
\right\}
\end{equation}

and $\{\zb\}$ is the vector containing the vertical coordinates of all the collocation points
(or degrees of freedom) of the BEM problem.  To better illustrate how the residual of the
numerical problem is put together, we split the vector in several parts. The vector $\hat{\phib}$
containing the values of the perturbation potential at the collocation nodes is split into
its portions $\hat{\phib}_D$ and $\hat{\phib}_N$ corresponding to points where Dirichlet or
Neumann boundary conditions are applied, respectively. The same kind of division is applied
to vector $\hat{\gammab}$, which its split into its degrees of freedom $\hat{\gammab}_D$ upon which
Dirichlet boundary conditions are applied, and its degrees of freedom $\hat{\gammab}_D$ upon which
Neumann boundary conditions are applied. Finally, the vector containing the vertical coordinates of
the collocation points $\{\zb\}$ is divided into its part containing the vertical coordinate of the
free surface and rest of the basin points $\zb_{FS}$ and $\zb_{BA}$, respectively. The
nonlinear system residual $\Fb$ is then split in a corresponding way into its portions $\Fb_{\phi_D}$,
$\Fb_{\phi_N}$, $\Fb_{\gamma_D}$, $\Fb_{\gamma_N}$, $\Fb_{z_{FS}}$ and $\Fb_{z_{BA}}$. Then
the residual components in each of these vectors are assembled in the following way

\begin{subequations}\label{eq:DAE_assembling}
\begin{eqnarray}
F_{\phi_D\ i} &=& \hat{\phib}_{D\ i} \\
F_{\phi_N\ i} &=& \alpha_i\hat{\phi_i} + \sum_{j=0}^{N_v} N_{ij} \hat{\phi}_j - \sum_{j=0}^{N_v} D_{ij} \hat{\gamma}_j \\
F_{\gamma_D\ i} &=& \alpha_i\hat{\phi_i} + \sum_{j=0}^{N_v} N_{ij} \hat{\phi}_j - \sum_{j=0}^{N_v} D_{ij} \hat{\gamma}_j \\
F_{\gamma_N\ i} &=& \sum_{j=0}^{N_v} M_{ij}\hat{\gamma}_j - b^{\phi_n}_i \\
F_{z_{FS}} &=& \sum_{j=0}^{N_v} \widetilde{M}_{ij}\hat{\dot{\phib}}_j - b^{\dot{\phi}}_i \\
F_{z_{BA}} &=& z_i - z_{REF\ i},
\end{eqnarray}
\end{subequations}
in which we made use of the vector $\zb_{REF}$ containing the value of the vertical coordinate of all the
nodes at the initial time step.

As will be discussed in more detail in
Section \ref{sec:test_case}, to test the effectiveness of the free surface boundary
condition formulation proposed and avoid other sources of error, in this work we
only consider the case of a fully submerged body advancing steadily in calm water. For such a reason
System~\eqref{eq:DAE} does not include any part that refer to the horizontal coordinates of
the collocation points, as they must experience no horizontal motion. 
As for the vertical coordinates,  System \eqref{eq:DAE_assembling} suggests that the
displacements --- and velocities --- will be set to
0 for all the collocation points, except for the ones on the free surface. In correspondence
with such nodes, the system equations will be obtained from System~\ref{eq:L2_proj_SUPG},
which represents the discretized and stabilized version of the ALE free surface
dynamic boundary condition. Note that the grid velocity field $\vb$ appearing in the 
ALE free surface and non homogeneous Neumann boundary conditions, is simply the time
derivative $\{\dot{\xb}\}$of the collocation point coordinates. Finally, as for the
the system degrees of freedom associated with
the potential $\hat{\phi}$ and potential normal derivative $\hat{\gamma}$ collocation point values,
the BEM resolution Equations~\eqref{LaplaceLinearSystem} are used.

Equation~\eqref{eq:DAE} represents a system of nonlinear differential algebraic
equations (DAE), which we solve using the IDA package of the
SUNDIALS OpenSource library \cite{sundials2005}. A $10^{-5}$ relative residual tolerance is
set for the Newton iterations used to solve the nonlinear problem
arising at each time step from the implicit time discretization scheme. In the linear
step of such iterations, the exact Jacobian of the numerical residual 
defined in Equation \eqref{eq:DAE_assembling} is considered. Such Jacobian, obtained
by means of automatic differentiation tools included in the package Sacado of
the C++ library Trilinos~\cite{Heroux2005} is inverted by means of a direct LU
factorization method. We point out that in the framework of IDA, once the solution at one time
step has been obtained, the initial guess at the ensuing time step is obtained
making use of the same BDF used in the implicit time advancing scheme. Typically, a
null solution is used at the initial time step.

{

We finally remark that in the framework of the DAE algorith employed, the ALE
time derivative of the velocity potential on the BEM collocation
points $\hat{\dot{\phib}}=\frac{\delta \hat{\phib}}{\delta t}$
is available not only at each time step, but also at each Newton correction.
Such ALE derivative yields the Eulerian derivative (see Equation~\eqref{eq:total_der})
at each collocation point. The velocity potential Eulerian time derivative is then plugged
into Bernoulli’s Equation~\eqref{eq:bern} to evaluate the pressure on the
whole domain boundary, without requiring the solution of
additional boundary value problems for $\frac{\Pa \phi}{\Pa t}$. The resulting
pressure ﬁeld can be integrated on the surface of any body of interest to
obtain the pressure force acting on it. We refer the interested reader
to \cite{molaIsope2016}, in which this approach was used to carry out
full fluid-structure interaction simulations of a ship free to move under the action
of hydrodynamic forces.
}

\subsection{DAE restart procedure}\label{sec:restarts}

A noteworthy feature of the solver developed, is that the time integration
is periodically paused to allow for adaptive grid refinement. At each
refinement cycle, a Kelly error estimator (\cite{kelly1983posteriori,gago1983posteriori,AinsworthOden-1997-a}) is
computed based on the water elevation field $\eta$. After the cells are
sorted according to the error indicator, a prescribed fraction of them
having the highest values are flagged and eventually refined.

Once the grid refinement has been carried out, all the fluid dynamic fields
are interpolated onto the new mesh. Of course, the interpolated solutions
will not satisfy the DAE residual in Equation~\eqref{eq:DAE}. Since
a non null initial residual normally leads to simulation blowup, the solution must
be adjusted at each restart, so as to satisfy Equation~\eqref{eq:DAE}.
In the model discussed in \cite{waveBem}, at each start, the solution $\overline{y}$
was obtained through interpolation of the coarse grid solution on the new grid. Then
the restart solution time derivative $\dot{\yb}^r$ was computed as the solution of the following nonlinear
equation system

\begin{equation}\label{eq:DAE_restart_vel}
F(\dot{\yb}^r,\overline{\yb},t^r)=G(\dot{\yb}^r)=0.
\end{equation}

System \eqref{eq:DAE_restart_vel} is solved by means of a Newton--Raphson algorithm implemented in the
KINSOL package of the SUNDIALS OpenSource library \cite{sundials2005}. Also in this case,
a $10^{-5}$ relative residual tolerance is set for the Newton iterations. And, also in this case,
the exact Jacobian of the numerical residual is computed, and inverted by means of a direct LU
factorization method.

The approach just described, which consists in imposing the interpolated nodes positions to obtain
the nodes velocities satisfying the DAE residual, did not lead to optimal results. In fact, it resulted
in very high nodes velocities which had to compensate for the slightly incorrect positioning of 
the nodes due to interpolation error. As a result, the time steps at each restart had to drop to follow
the faster dynamics, slowing down the simulation.

A much more interesting alternative restart treatment, is that of imposing the interpolated solution time
derivative $\overline{\dot{\yb}}$, to obtain the solution from the DAE residual, namely

\begin{equation}\label{eq:DAE_restart_pos}
F(\overline{\dot{\yb}},\yb^r,t^r)=G(\yb^r)=0.
\end{equation}

This approach does not introduce spurious fast dynamic components, and is of course to be preferred.
Unfortunately, numerical evidence suggests that if used --- as is the case for \cite{waveBem} --- in presence of the
semi-Lagrangian kinematic and dynamic boundary conditions \eqref{eq:kin_fs_condition_semil}
and \eqref{eq:dyn_fs_condition_semil}, Problem \eqref{eq:DAE_restart_pos} is not well posed,
likely due to a singular Jacobian $\frac{\Pa F}{\Pa \yb}$. On the other hand, the introduction
of the alternative formulation adopted in the present work, allows for the solution of
Problem \eqref{eq:DAE_restart_pos}, obtaining correct restart solution without introducing
spurious faster dynamics into the DAE system. In addition, as will be explained in next
session, the correct solution of Problem \eqref{eq:DAE_restart_pos} has been the gateway
to the possibility of solving steady problems.

\subsection{Stationary solver}\label{sec:steady}

As discussed at length, both the mathematical formulation
of the free surface boundary condition, and the numerical discretization of the
resulting boundary value problem have been selected so that stationary and
non stationary problems could have a unified implementation. 
To write the problem for the stationary solution $\yb^s$, we 
introduce the additional conditions $t\rightarrow\infty$ and $\dot{\yb}=0$
in System~\eqref{eq:DAE}, namely

\begin{equation}\label{eq:DAE_stationary}
F(\cancelto{0}{\dot{\yb}},\yb^s,t\rightarrow\infty)=G(\yb^s)=0.
\end{equation}
To obtain a unified implementation for steady and unsteady solvers, the nonlinear System~\eqref{eq:DAE_stationary} is solved
making use of the same residual function implemented for the DAE solver, in which the
$\dot{\yb}$ argument is set to zero at every call. The resulting problem is by all means
a particular case of Problem \eqref{eq:DAE_restart_pos}, and such a nonlinear 
system of equations for $\yb^s$ is again solved by means of a Newton--Raphson algorithm implemented in the
KINSOL package of the SUNDIALS OpenSource library. As is the case for the nonlinear problem associated with
restarts, the initial guess of the Newton iterations for the steady problem is obtained from the
previous adaptive step solution interpolated onto the new grid.
As will be shown, this resulted in a software in which it is possible to switch from
non stationary to stationary solver at the sole cost of including or not including time
derivatives upon numerical resolution of the DAE system. This is of course more straightforward than using
a set of completely different non linear free surface boundary conditions for unsteady and steady
potential flow problems. A flow chart representing the numerical procedure used in the software
is presented in Figure \ref{fig:flow_chart}.

\section {Results}\label{sec::results}

A simulation campaign has been carried out to fully characterize the performance
of the algorithm proposed. In particular, to reproduce the possible practical use of
a potential flow solver suited for early design stages, all the numerical tests have been
carried out on an Intel Quad Core i7-7700HQ 2.80GHz, 32 GB RAM laptop using 10 parallel processors.
In addition, the academic test case considered allowed for cross validation through the comparison with
well assessed literature results. The next sections will describe the details of the
test case considered, and present the results of the simulation campaign.

\subsection {The immersed ellipsoid test case description} \label{sec:test_case}

The test case considered is that of a fully immersed ellipsoid advancing at steady
speed in calm water. The spheroid considered is moving in the direction of its
horizontal axis of revolution and has a radius which is one fifth of its length.
Figure~\ref{fig:domain_dimensions} displays a two dimensional diagram reporting the
dimensions of the computational domain employed throughout the simulation campaign,
which is attached to the spheroid and is advancing in the water alongside with it. 
Most of the lengths reported in the picture are referred to the ellipsoid length  $L=10\ $m.
We also report that the overall width of the channel, which is not appreciable from this two
dimensional sketch, was set to $10L$.

Finally, the distance from the origin at which the
damping zone starts acting and the overall length of the damping zone have been set to
$x_d=50\,$m and $L_d=100\,$m.

\begin{center}
\begin{figure}
\centerline{
  \ifpdf
  \resizebox{1\textwidth}{!}{
    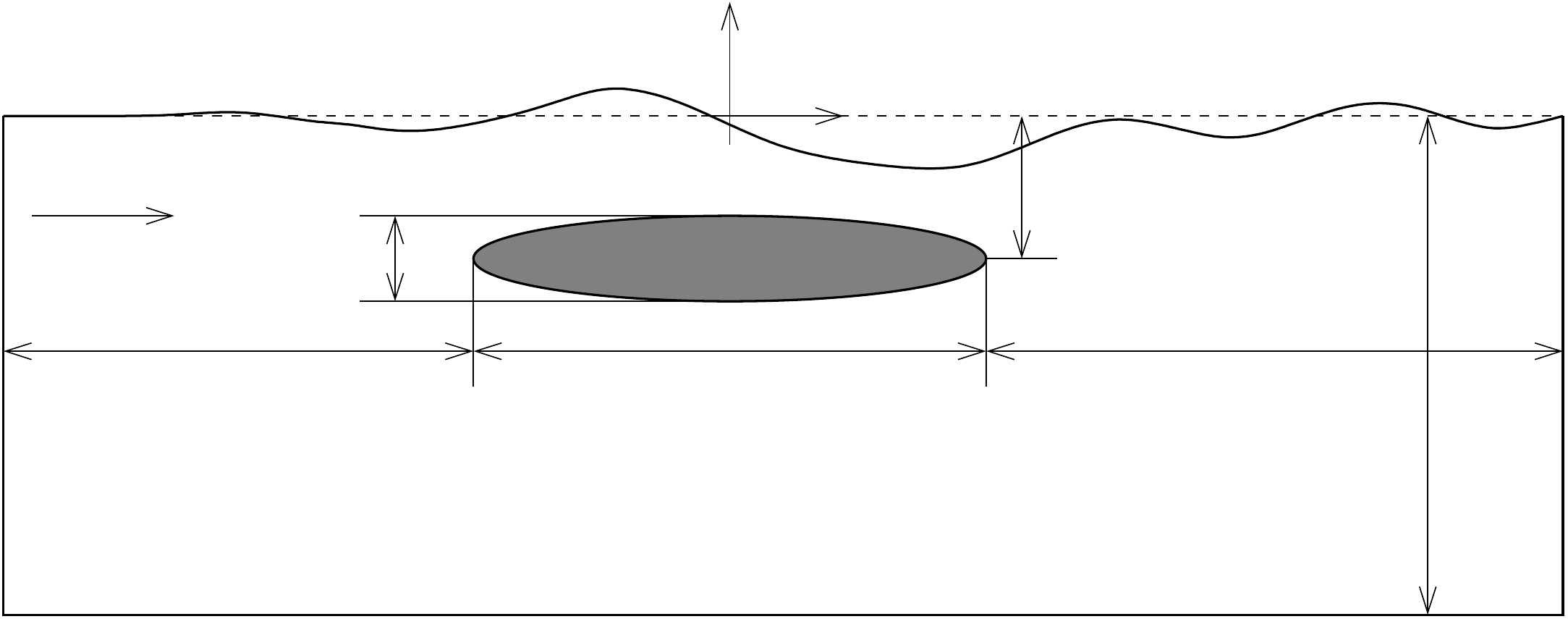
  }
  \else
  \resizebox{1\textwidth}{!}{
    \input{./figures/spheroid_domain.pstex_t}
  }
  \fi
}
\caption{A two dimensional sketch illustrating the computational domain dimensions used for the fully
          immersed ellipsoid test case. All the lengths indicated in the diagram are referred to
         the ellipsoid length $L=10\ $m. In the simulation campaign, we have tested different
         values of the vertical distance $f$ between the spheroid horizontal axis of revolution and
         the undisturbed free surface elevation, located at $z=0$.\label{fig:domain_dimensions}}
\end{figure}
\end{center}

\subsubsection{Steady state numerical experiments}

A first set of experiments has been devoted to evaluate the performance of the
steady flow solver developed, both in terms of computational cost and in terms of
accuracy of the solutions. In this case, the asymptotic potential is set to $\phi_\infty = \Ub_\infty\cdot\xb$,
and a constant velocity value $\nablab \phi_\infty =\overline{\Ub}_\infty$ is imposed in each simulation considered.
In such set of tests, several values of $\overline{\Ub}_\infty$ aligned with the $x$ axis of the domain
have been considered, as well as
multiple values of the depth $f$ --- defined as the vertical distance between the spheroid
horizontal axis of revolution and the undisturbed free surface elevation, located at $z=0$.
For a better evaluation of the results, the non dimensional version of the latter parameters
will be reported in the next sections. The Froude number (namely $\text{Fr}=|\Ub_\infty|/(g L)$ ) will be
used as the non dimensional measure of the asymptotic velocity, while the non dimensional ellipsoid
distance from the undisturbed free surface will be indicated by the parameter $d/f$.

Making use of the CAD handling features of the $\pi$-BEM library \cite{pi-BEM}, the very
coarse quadrilateral mesh originally imported is automatically refined on the surface spheroid
until it suitably represents the object geometry. 

\begin{center}
  \begin{figure}
  \includegraphics[width=0.48\linewidth]{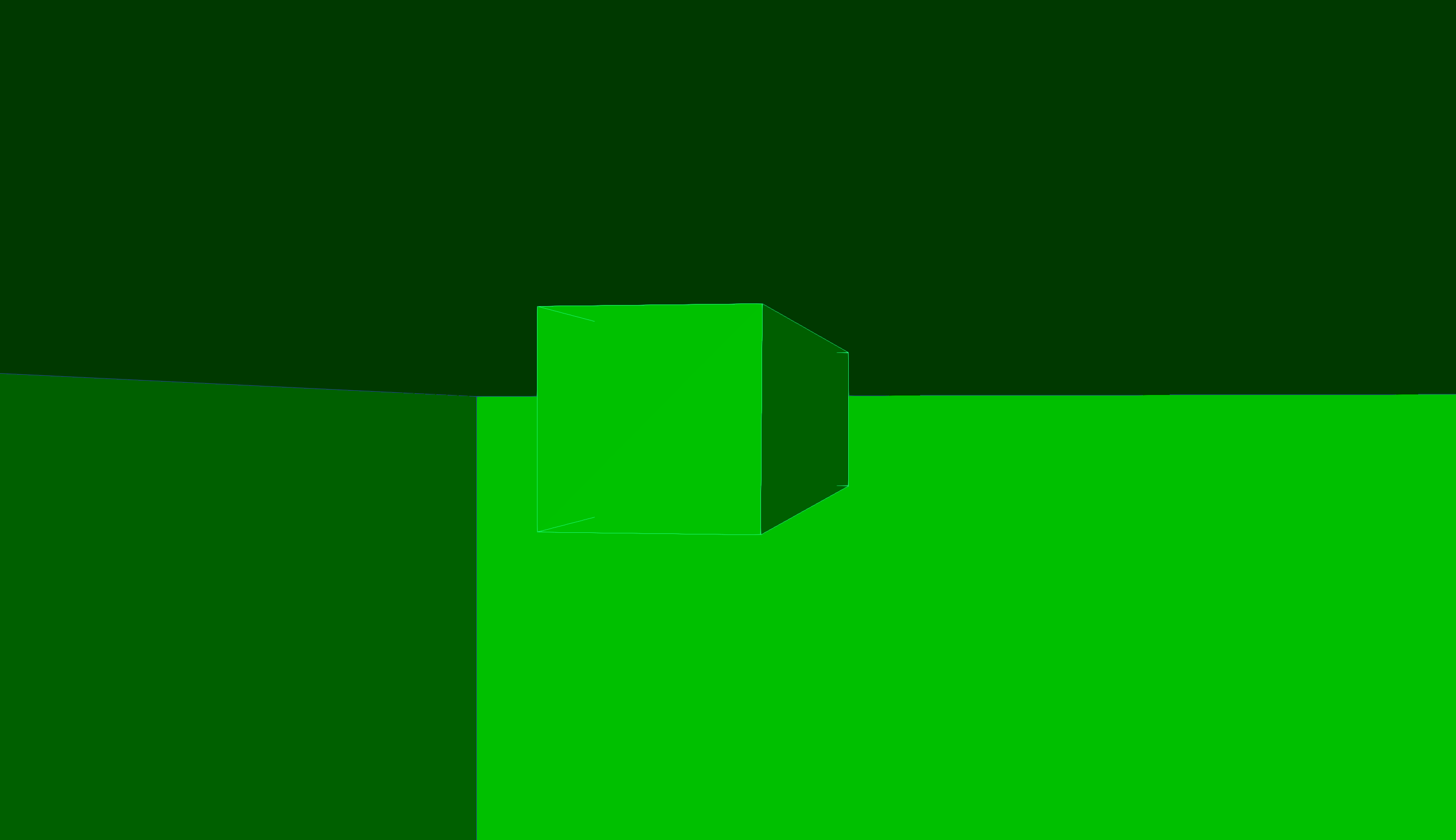}
  \includegraphics[width=0.48\linewidth]{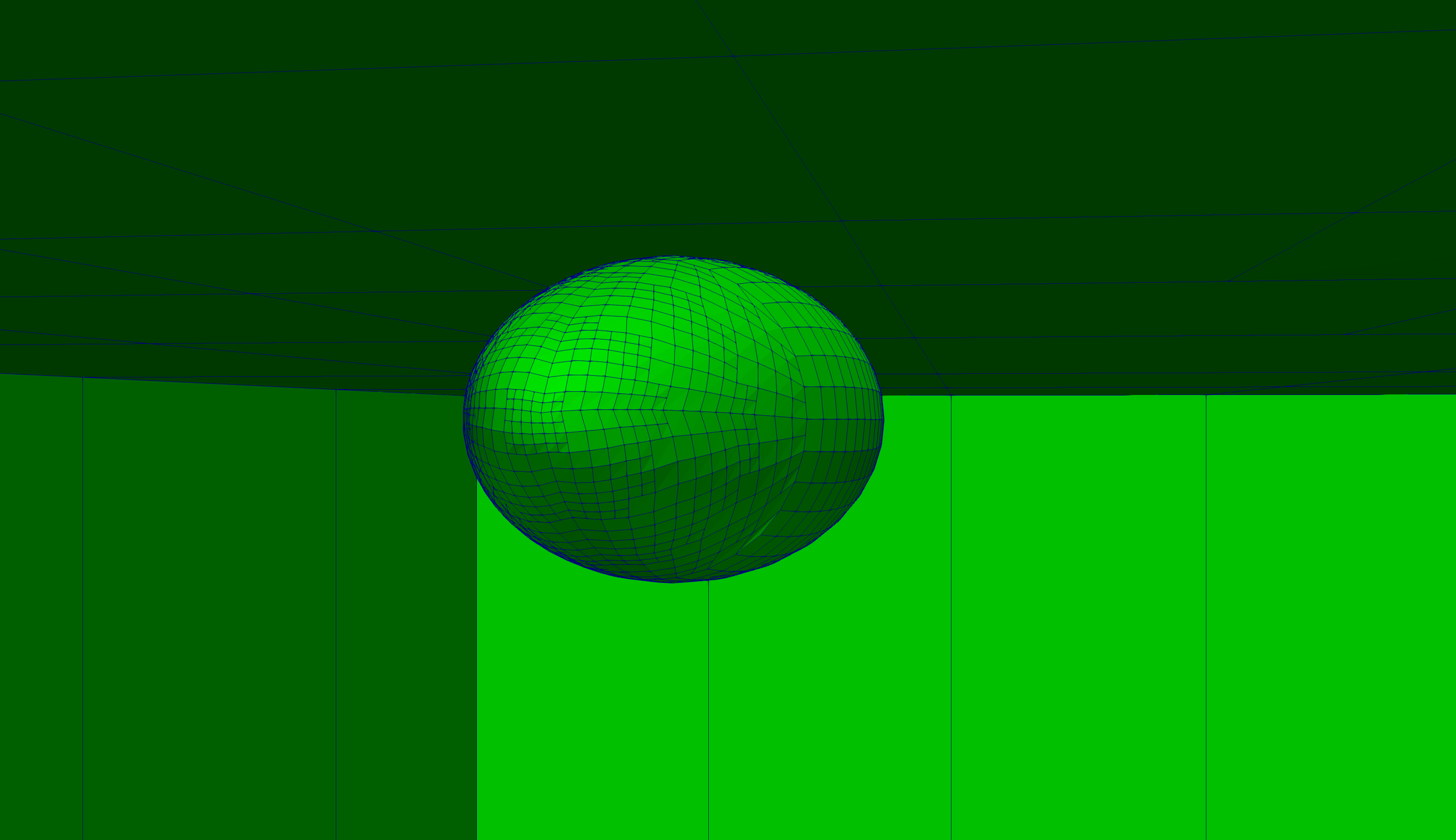}
  \caption{On the left, a view of the ellipsoid initial mesh, only composed of 4 cells. On the
           right, the same view of the ellipsoid mesh used for all the simulations. The latter
           mesh has been obtained through fully automated adaptive refinement cycles based on
           surface curvature, and interfaced with the ellipsoid CAD geometry.\label{fig:spheroid_mesh}}
  \end{figure}
  \end{center}

Figure~\ref{fig:spheroid_mesh} illustrates
such process, which is made up of 7 adaptive refinement cycles based on the curvature. Throughout
each refinement cycle, the CAD geometry is interrogated to compute the position of the
new nodes on the spheroid surface, and to obtain evaluation of local curvature to flag cells for
refinement. In such framework, the original coarse mesh (on the left in the figure)
composed of 16 nodes is refined to obtain the simulation initial mesh featuring 1989 nodes
(on the right in the figure). We also point out that the automated refinement process includes
cycles to reduce the cells aspect ratios until they are lower than 3.5, and cycles to
refine the region of the free surface closest to the ellipsoid.

The initial computational mesh described is depicted in Figure~\ref{fig:free_surface_mesh}, on the
left. The other plots in the Figure refer to further cycles of the adaptive refinement process
based on surface elevation. In such process, once the initial mesh is available, the nonlinear
problem resulting from the steady fluid dynamic equations is solved to compute the flow velocity
potential and the water elevation. The latter field is then used to compute Kelly error estimator
and flag for refinement the  top 4\% portion of free surface cells. The cycle is then repeated
 15 additional times, to obtain the  grids depicted
on the right of Figure~\ref{fig:free_surface_mesh}. As expected, by a qualitative standpoint
the computational grid refinement pattern appears to follow the V-shaped Kelvin wake induced
by the spheroid underwater motion. 

\begin{center}
\begin{figure}
\begin{tabular}{c c c c}
\hspace{-1.3cm} \includegraphics[width=0.39\linewidth]{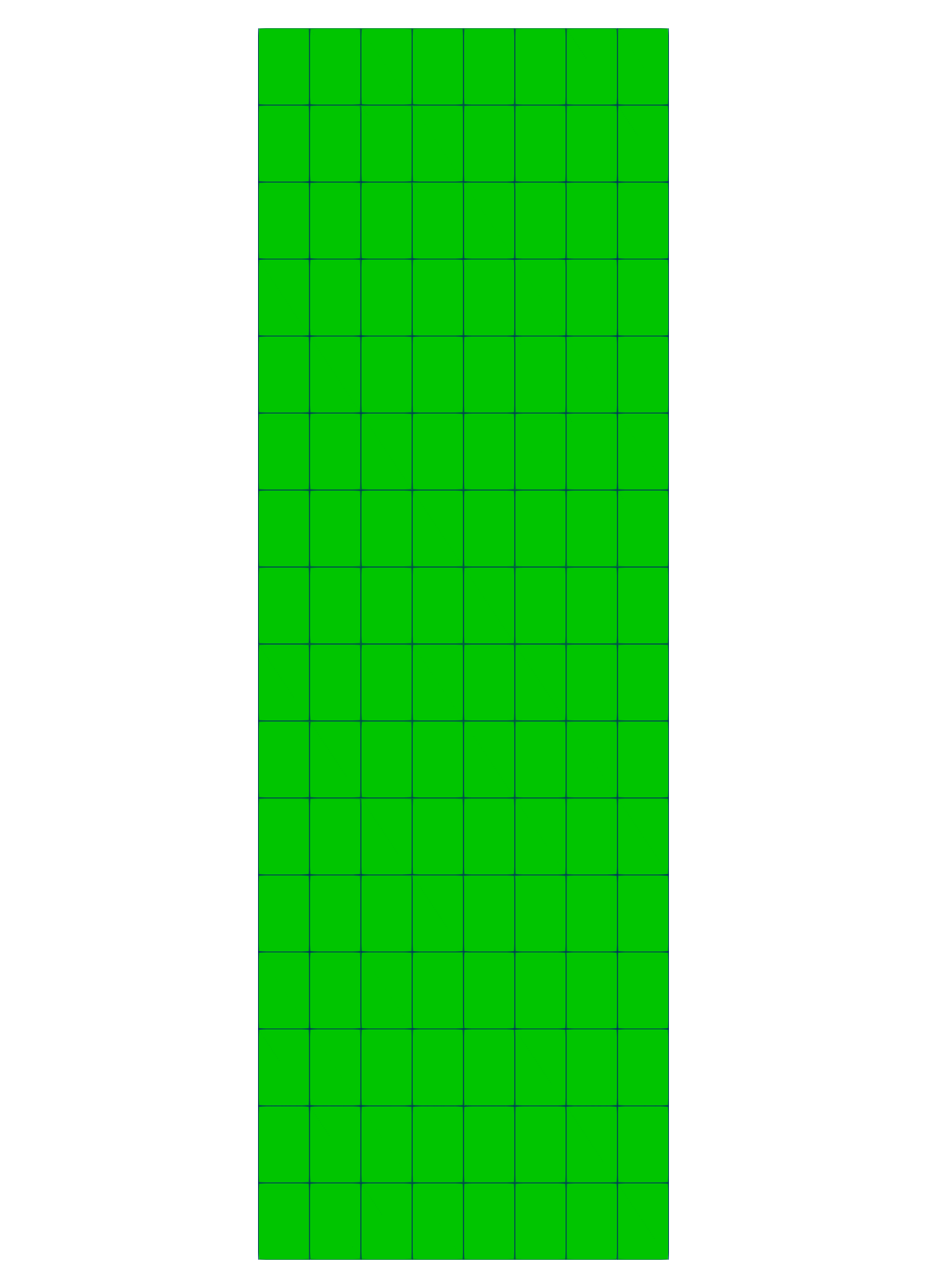} \hspace{-1.2cm} &
\hspace{-1.2cm} \includegraphics[width=0.39\linewidth]{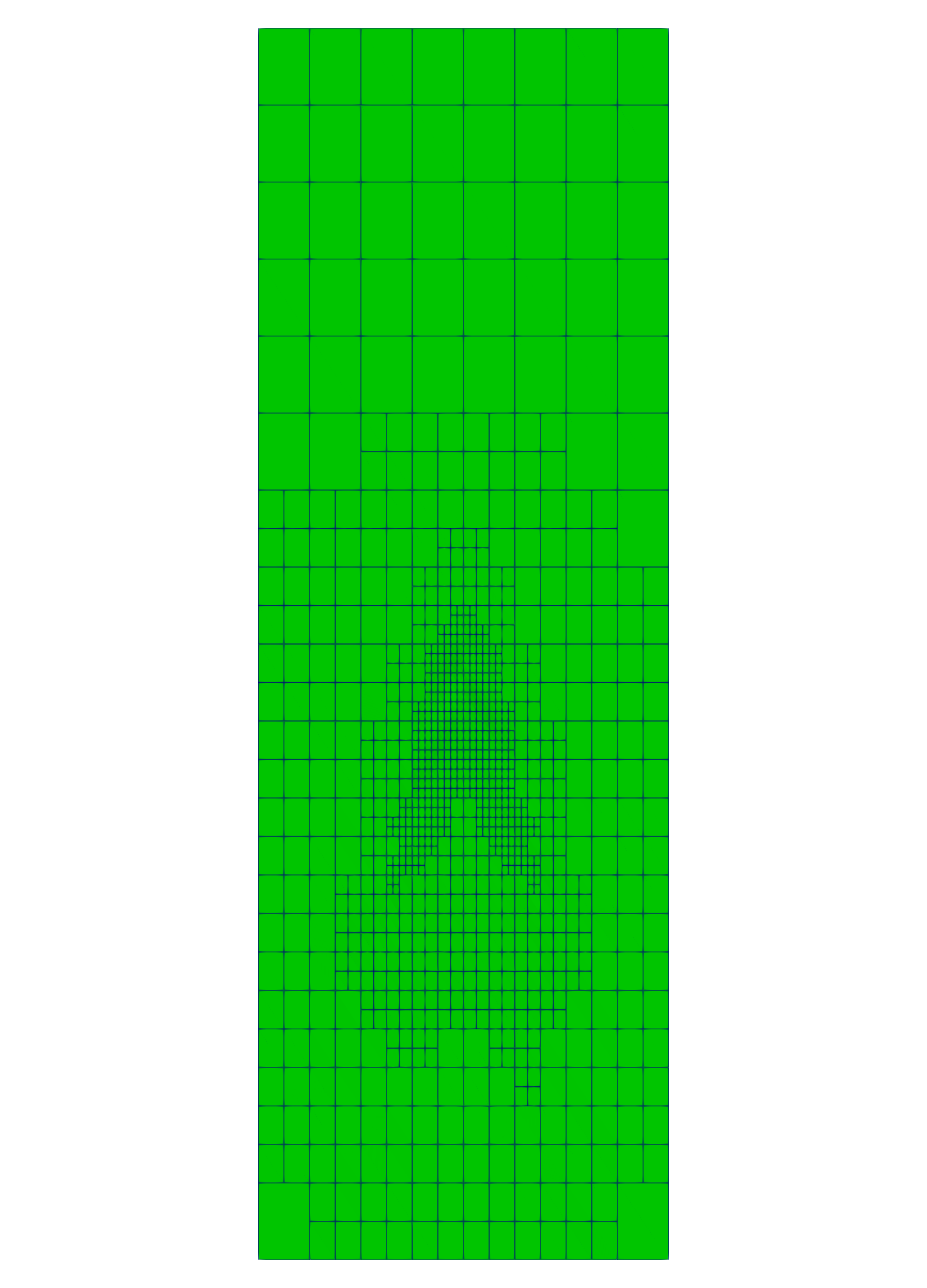} \hspace{-1.2cm} &
\hspace{-1.2cm} \includegraphics[width=0.39\linewidth]{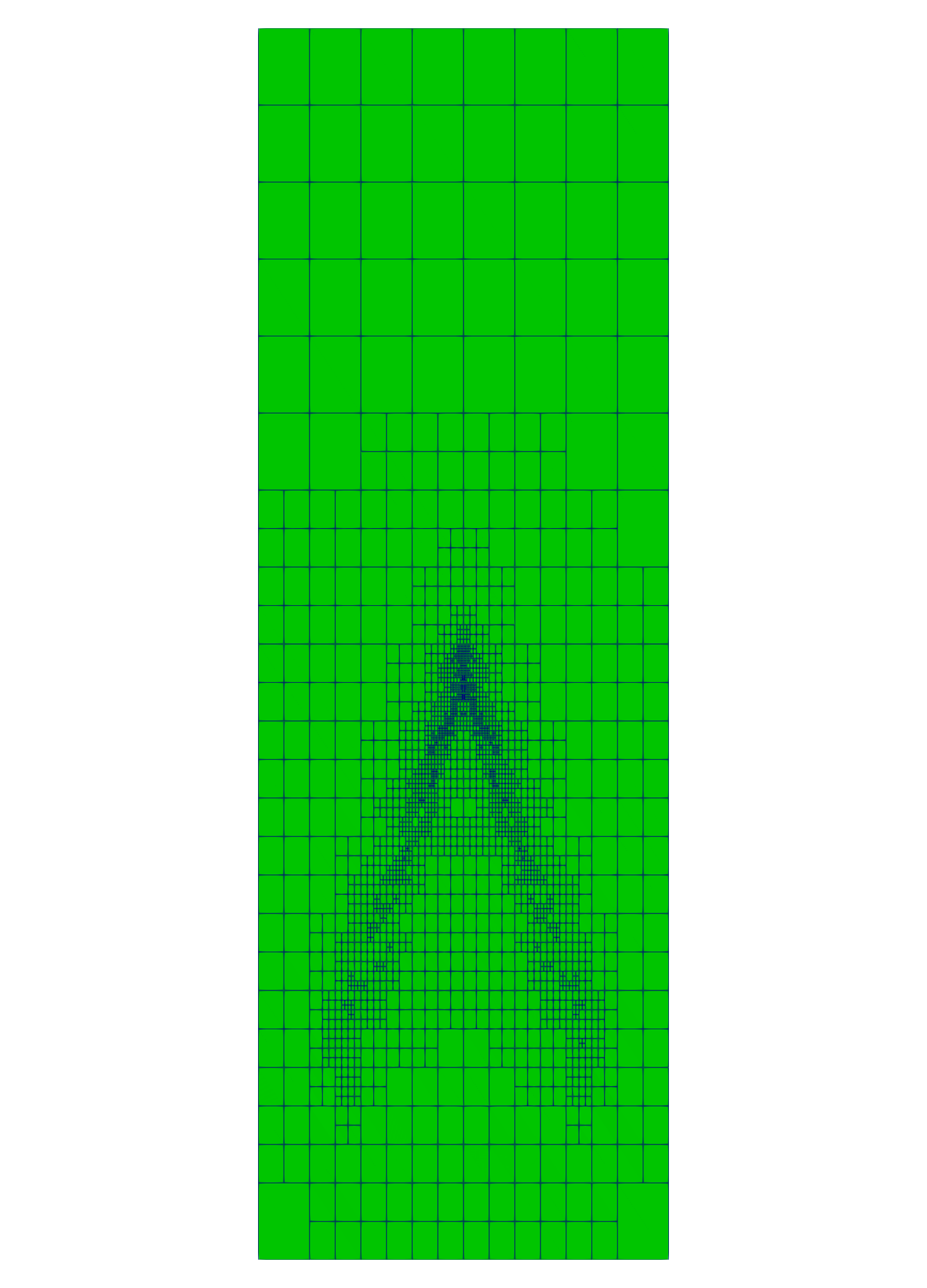} \hspace{-1.2cm} &
\hspace{-1.2cm} \includegraphics[width=0.39\linewidth]{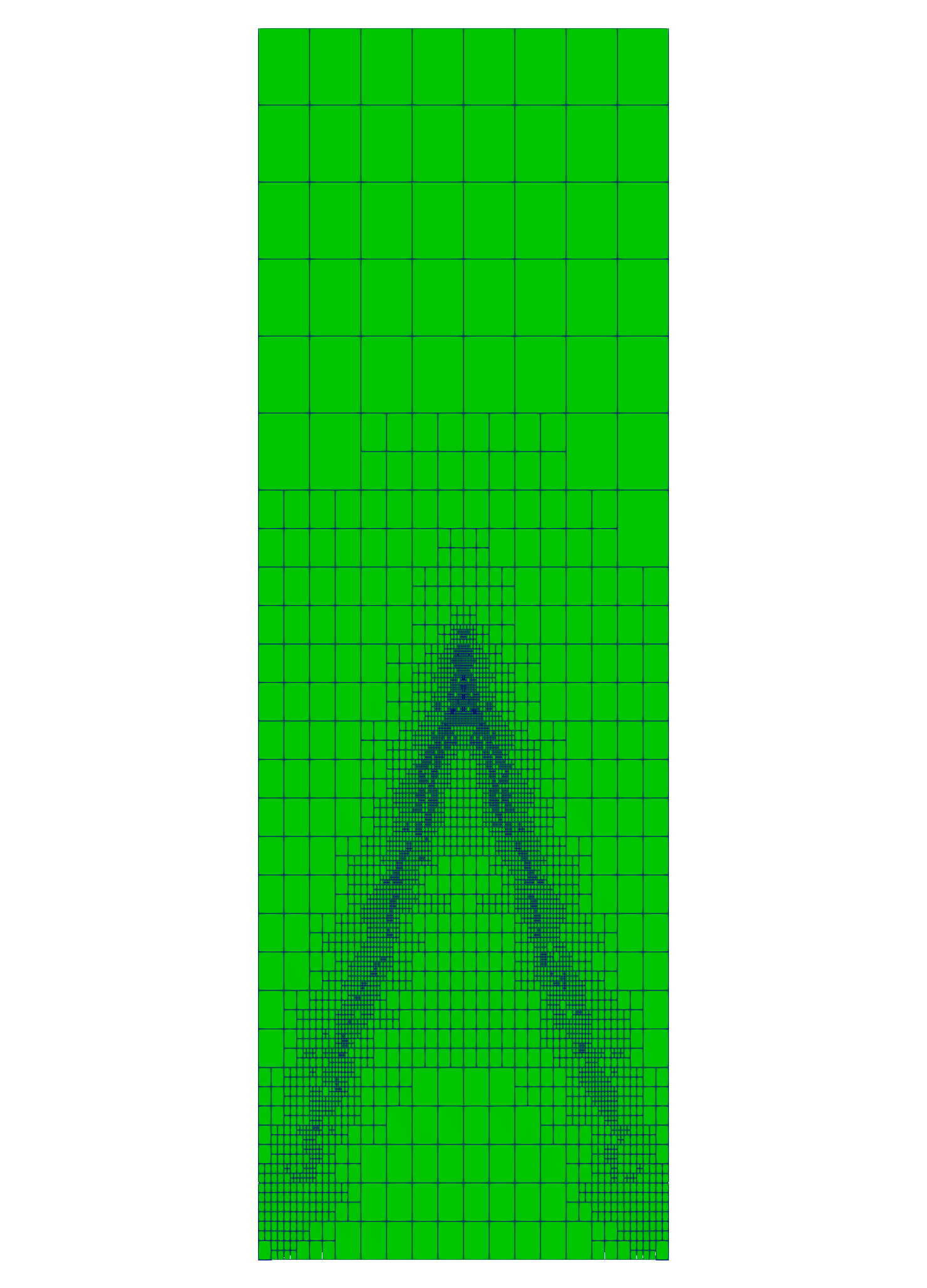} \hspace{-1.2cm}\\
Initial & Cycle 4 & Cycle 8 & \hspace{-1.3cm}Cycle 12
\end{tabular}
\caption{An illustration of the computational grids across the adaptive refinement process based
         on free surface elevation characterizing the simulations. From left to right, the top
         view of the free surface grid at Cycles 0, 4, 8, 12 respectively. As expected, the refinement
         pattern covers the Kelvin wave produced by the spheroid.\label{fig:free_surface_mesh}}
\end{figure}
\end{center}

\subsubsection{Numerical experiments of unsteady flows having steady regime solution}\label{sec:fake_unsteady_experiments}

The unsteady flow simulations carried out in this work have been aimed at providing a
numerical confirmation that the solver developed can indifferently operate under transient or
steady flow assumptions. In particular, it is relevant showing that if a flow admits
a regime solution that is constant over time, both the transient and the steady solver converge
to the same result. Thus, to allow for such a comparison, we considered one of the steady test
cases previously described, and solved with the transient solver. We focused in particular
on the case in which the depth of the cylinder corresponded to $d/f = 0.8$ and $\text{Fr}=0.8$.
In the numerical experiments, we considered three different transient simulations in which the
asymptotic velocity $\Ub_\infty$ aligned with the $x$ axis of the domain is increased with sinusoidal
ramps lasting $0.75\,s$, $7.5\,s$ and $15\,s$,
respectively.  The aim of the present test is that of providing numerical confirmation to the
fact that the same numerical residual is solved by both the steady and unsteady solver in the
limit $t\rightarrow\infty$, which results in the same steady and unsteady solver solution. Thus, to avoid
that different adaptive refinement patterns in the steady
and unsteady solver could introduce even a small error in such evaluation, a steady simulation
featuring 4 adaptive refinement cycles is first carried out. After this, the solution is reinitialized
on the adaptively refined grid, and a transient simulation is started from rest, and run until
a convergence to a regime solution is observed. This procedure allows
for the evaluation of the wave resistance and lift forces obtained both with the steady and unsteady
solver, when the same grid is used. Since, as explained in Section \ref{sec:steady}, the nonlinear
problem residual used by both methods is the same for $t\rightarrow \infty$, the two solutions are expected to be
identical for $t\rightarrow \infty$.

\subsubsection{Numerical experiments of unsteady flows having unsteady regime solution}\label{sec:true_unsteady_experiments}

A full investigation on the influence of direction, amplitude and wavelength of incident waves
on the wave resistance of a hull falls way beyond the scope of the present work. However,
to assess the ability of the present free surface model to simulate the flow past a ship hull advancing in a wave field,
we also consider completely unsteady test case. Here, the asymptotic potential $\nablab\phi_\infty$ is the sum of
an Airy linear wave potential and a constant velocity potential, namely

\begin{equation}\label{eq:unst_asym_pot}
\phi_\infty(\xb,t) = U_\infty x + a\frac{\omega}{k}\frac{\cosh(k(z+h))}{\sinh(kh)}\sin(kx-\omega t).
\end{equation}
Also in this case, both the Airy wave and the constant velocity vector $\Ub_\infty$ are aligned with the
$x$ axis of the computational domain. In \eqref{eq:unst_asym_pot} the velocity magnitude is set to
$U_\infty=2.9714\,$m/s corresponding to $Fr=0.3$, the incident wave amplitude is $a=0.02\, $m,
and the wave length is $\lambda=4L=40\,$m corresponding to $k=0.1571\,$1/m. Considering the height of
the basin $h=50\,$m, the dispersion relation $\omega^2=gk\tanh(kh)$ results in $\omega = 1.2414\,$rad/s
and a corresponding period of $T=5.0616\,$s. The asymptotic free surface elevation considered in this
numerical test is

\begin{equation}\label{eq:unst_asym_eta}
\eta_\infty = a\cos(kx-\omega t),
\end{equation}
which is identical to the Airy wave elevation. We must point out that in principle, making use of the linearized
kinematic free surface condition should result in additional cross terms to account for the presence of a
stream velocity. Such terms are not accounted for in this work, which for the case at hand results in variations
of the Bernoulli constant for the asymptotic flow $C(t)$ evaluated on the free surface of approximately $\pm 1\%$ with respect to
its average value. Once again, the main goal of the present investigation in only assessing whether a set
of waves with the desired wavelength, amplitude and period can be introduced in the numerical domain and can
interact with the Kelvin wave pattern generated by the hull. For such a reason, in the framework of the
present investigation we considered such an approximation acceptable. Further investigations will be addressed at
forcing wave fields with different desired characteristics in the numerical towing tank developed.

The numerical simulation of the unsteady test cases described
have been obtained making use of a mesh in which the vertical position of the ellipsoid with respect to the
undisturbed free surface corresponds to $d/f = 0.75$. In addition, a different set of initial refinements has been
considered, so as to be able to cover each wavelength of the incident waves with at least four cells. From such
grid, a steady test case featuring 7 refinement cycles was carried out to obtain a grid able to properly describe
a Kelvin wake pattern. The whole procedure resulted in a grid of approximately 4000 nodes. The simulation was
then restarted from a solution corresponding to null velocity and
and water elevation, and the asymptotic potential was then increased to its regime values by means of a $5\,$s
sinusoidal ramp. Also in this case, the procedure was devised so as to obtain the wave resistance for the steady
and unsteady case on the same grid, consequently avoiding possible differences associated with different grids. In
addition, avoiding adaptive refinements diring the time integration of the unsteady problem limited the computational
cost.

\subsection{Numerical results}

A typical output of the simulations
is portrayed in Figure~\ref{fig:free_surface_soluton}. The picture refers to
the $d/f = 0.8$ and $\text{Fr}=0.6$ test case, and both the free surface
and underwater spheroid mesh are visible. The free surface is colored
according to contours of water elevation, which make visible the accurate
reconstruction of the Kelvin wake pattern generated by the moving body.

\begin{center}
\begin{figure}
\includegraphics[width=0.9\linewidth]{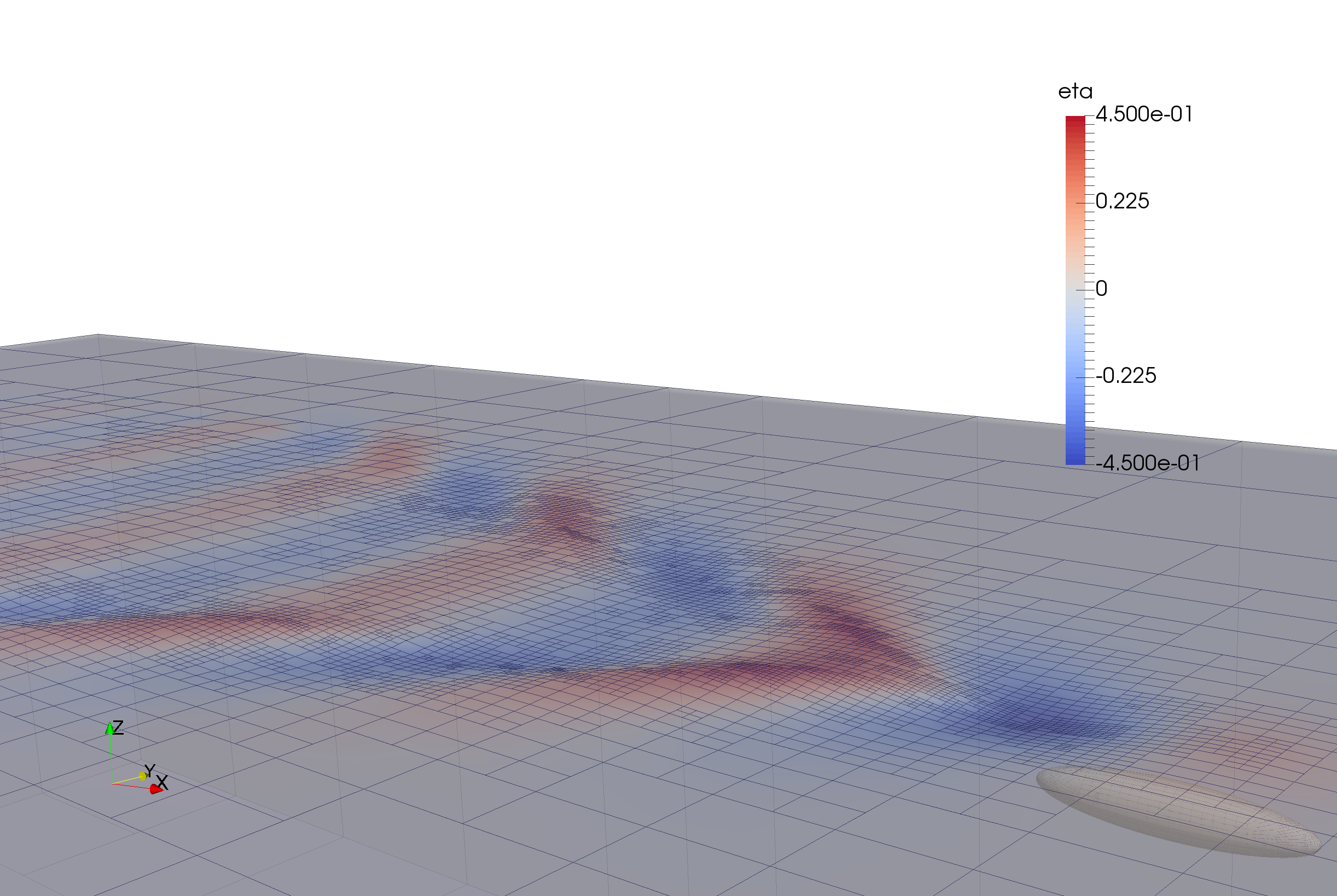}
\caption{An illustration of the free surface elevation contour, and of the domain mesh. The plot refers
         to the $d/f = 0.8$ and $\text{Fr}=0.6$ test case.\label{fig:free_surface_soluton}}
\end{figure}
\end{center}

We will now present the quantitative results of the campaign of numerical
experiments described in Section~\ref{sec:test_case}. A first objective will
be that of characterizing the performance of the adaptive refinement approach
adopted in terms of both computational cost and convergence to a stable solution.
After this, we will discuss the accuracy of the solutions obtained.

\subsubsection{Steady case results}

The overall computational cost of a typical simulation is reported in
Table \ref{tab:time_report} --- which more specifically refers to the stationary
ellipsoid test case in which $\text{Fr}=0.7$ and $d/f = 0.8$.
For each line corresponding to a refinement
cycle carried out, the columns of the Table report the number of computational
grid nodes, along with the number of nonlinear iterations and Jacobians assemblies
required by the Newton solver to reach convergence, and finally the execution time.
As expected, the dimension of the grid grows at an approximately constant rate in at each refinement
cycle. This is explained by the fact that only the portion of cells with highest
error indicators is
refined at every cycle. Thus, as more cycles are executed, such error gets distributed
over a wider amount of cells, which will be then flagged for refinement. For the most
part, the  execution time of each cycle reported in the table mirrors the growth in the
computational grid nodes --- which are also the BEM problem degrees of freedom. The only
factor perturbing the expected quadratic relationship between grid nodes and computational
cost, is the number of nonlinear iterations and of Jacobians assemblies
required by the Newton solver to converge at every cycle. The latter operation is in fact
particularly time consuming, and can significantly affect the duration of a single simulation
cycle. Typically, the Jacobian is assembled once per each refinement cycle, but because the cost
of its computation and LU factorization grows as the third power of the number of degrees of freedom of
the problem, its impact grows at each refinement cycle. The overall
computational cost of the entire refinement cycles procedure is approximately five hours,
which is in principle compatible with simulations run in early hull design stages. We must also
remark that attempting additional refinement cycles past the 12th one results in an arrest of
the computation. Rather than being associated with a divergence of Newton iterations, such an
arrest is due to a failure of the Jacobian matrix LU factorization. This is consistently
observed for computational grids with more than 9000 nodes, which result in Jacobian systems
of more than 36000 unknowns. It is possible at such size, the condition number of
the --- for the most part dense --- Jacobian matrix could become too big or too ill conditioned
for the factorization algorithm to properly work. Future work will be devoted to study appropriate preconditioning
strategy of the Jacobian matrix that would allow for the use of Krylov subspace based linear solvers,
more suited for the dimension of the linear system.

\begin{table}
\begin{center}

\begin{tabular}{||c|c|c|c|c||}
\hline
\hline
Refinement       & Grid       & Nonlinear           & Jacobians           & Cumulative                               \\
cycle            & nodes      & iteratios           & assembled           & execution time                                    \\
\hline
0                &  1989      &  4                  &  1                  &  0\,\text{h} 01\,\text{m} 38\,\text{s}  \\
\hline
1                &  2260      &  6                  &  1                  &  0\,\text{h} 03\,\text{m} 47\,\text{s}  \\
\hline
2                &  2568      &  6                  &  1                  &  0\,\text{h} 06\,\text{m} 53\,\text{s}  \\
\hline
3                &  2865      &  6                  &  1                  &  0\,\text{h} 11\,\text{m} 10\,\text{s}  \\
\hline
4                &  3207      &  5                  &  1                  &  0\,\text{h} 16\,\text{m} 55\,\text{s}  \\
\hline
5                &  3621      &  5                  &  1                  &  0\,\text{h} 24\,\text{m} 49\,\text{s}  \\
\hline
6                &  4086      &  6                  &  1                  &  0\,\text{h} 36\,\text{m} 00\,\text{s}  \\
\hline
7                &  4621      &  5                  &  1                  &  0\,\text{h} 51\,\text{m} 36\,\text{s}  \\
\hline
8                &  5213      &  5                  &  1                  &  1\,\text{h} 12\,\text{m} 56\,\text{s}  \\
\hline
9                &  5896      &  5                  &  1                  &  1\,\text{h} 43\,\text{m} 08\,\text{s}  \\
\hline
10               &  6642      &  5                  &  1                  &  2\,\text{h} 25\,\text{m} 37\,\text{s}  \\
\hline 
11               &  7529      &  5                  &  1                  &  3\,\text{h} 26\,\text{m} 16\,\text{s}  \\
\hline
12               &  8477      &  5                  &  1                  &  4\,\text{h} 52\,\text{m} 00\,\text{s}  \\
\hline
\hline

\end{tabular}

\end{center}
\caption{Grid nodes, number of Jacobian matrices assembled and overall computational times required
         for each adaptive refinement cycles. The values reported refer to the stationary
         ellipsoid test case in which $\text{Fr}=0.7$ and $d/f = 0.8$, solved on an Intel Quad Core i7-7700HQ 2.80GHz, 32 GB RAM laptop
         using 10 parallel processors. \label{tab:time_report}}
\end{table}

Figure~\ref{fig:ForceConvergence} allows for an evaluation of the adaptive refinement cycles
effectiveness in converging to a stable solution. The left diagram in the Figure shows the
typical evolution of the steady state hydrodynamic lift across 12 adaptive refinement cycles.
The plot refers to the stationary test case in which $\text{Fr}=0.7$ and $d/f = 0.8$. As can be
appreciated, the hull hydrodynamic lift gradually decreases to values lower than the
hydrostatic lift ($L_0 = 209010.89\ $N), and appears to settle in the last two iterations
to values that are approximately 1.3\% shorter than $L_0$.

\begin{center}
\begin{figure}
\includegraphics[width=0.49\linewidth]{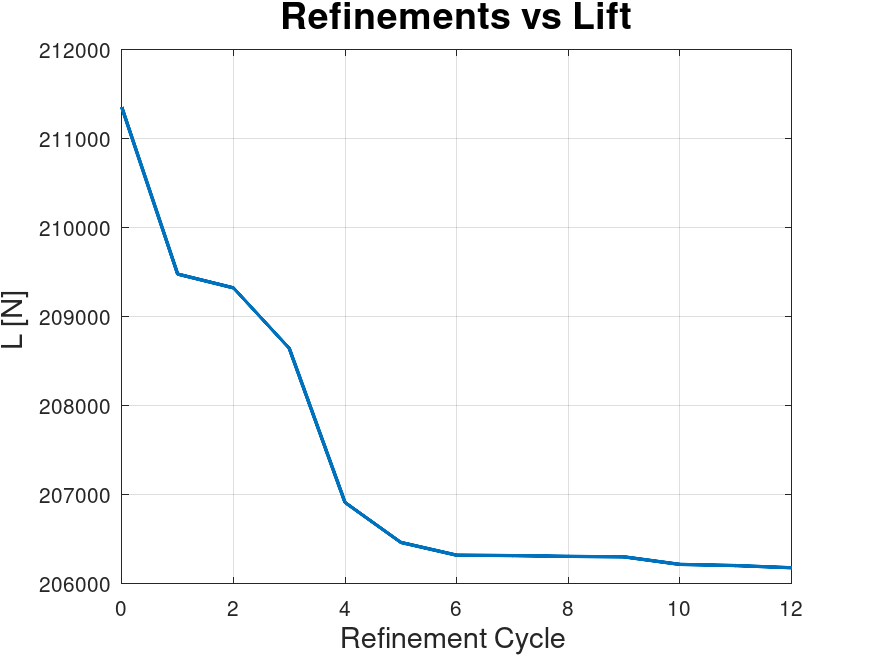}
\includegraphics[width=0.49\linewidth]{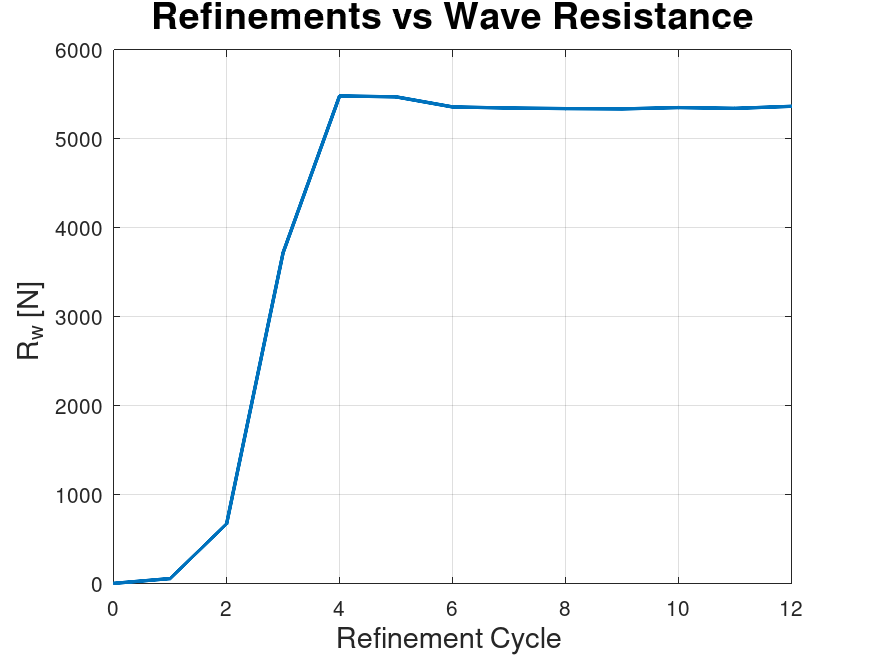}
\caption{Fully immersed ellipsoid dimensional lift (left) and resistance (right)value as
         a function of adaptive refinement iterations. The curve refers to the stationary
         ellipsoid test case in which $\text{Fr}=0.7$ and $d/f = 0.8$.\label{fig:ForceConvergence}}
\end{figure}
\end{center}

The grid convergence trend is further confirmed by the  right plot in Figure~\ref{fig:ForceConvergence},
which depicts the wave resistance evolution across 12 refinement cycles. The plot refers again to the stationary
test case in which $\text{Fr}=0.7$ and $d/f = 0.8$. Starting from low values obtained with the initial coarse
grids, the resistance values gradually increase across refinements, until the last refinement cycle
considered result in no significant resistance variation.

The results presented in
Figure~\ref{fig:ForceConvergence} suggest that for the ellipsoid test
case under study grid convergence is obtained at the 6th refinement cycle. In fact, the wave
resistance and lift forces computed at the 6th refinement cycle only differ by  0.13\% and
0.07\%, respectively, from the values computed at the 12th and final refinement. Based on the
computational costs reported in Table \ref{tab:time_report}, it can be inferred that a reliable
drag and lift prediction is obtained in approximately 36\,\text{m} on a laptop.

Figure~\ref{fig:ForceCompareAcrossRefs} provides a
confirmation that such trend holds across all the range of Froude numbers investigated.
The plots in the Figure show the values on non dimensional net lift $L^*=(L-L_0)/L_0$ (left)
and non dimensional wave resistance $R^*=(R)/L_0$ (right) as a function of Fr. In both diagrams
the diamonds indicate the results obtained at the 6th adaptive refinement cycle, while the
asterisks refer to the results obtained at the  12th cycle. For reference, continuous lines
representing corresponding literature results by Scullen~\cite{Scullen} have been added to the
plot. As can be appreciated in all the test cases considered the difference between the
solution at the last two refinement cycles is minimal, even compared to the difference
observed with different models solutions.

\begin{center}
\begin{figure}
\includegraphics[width=0.49\linewidth]{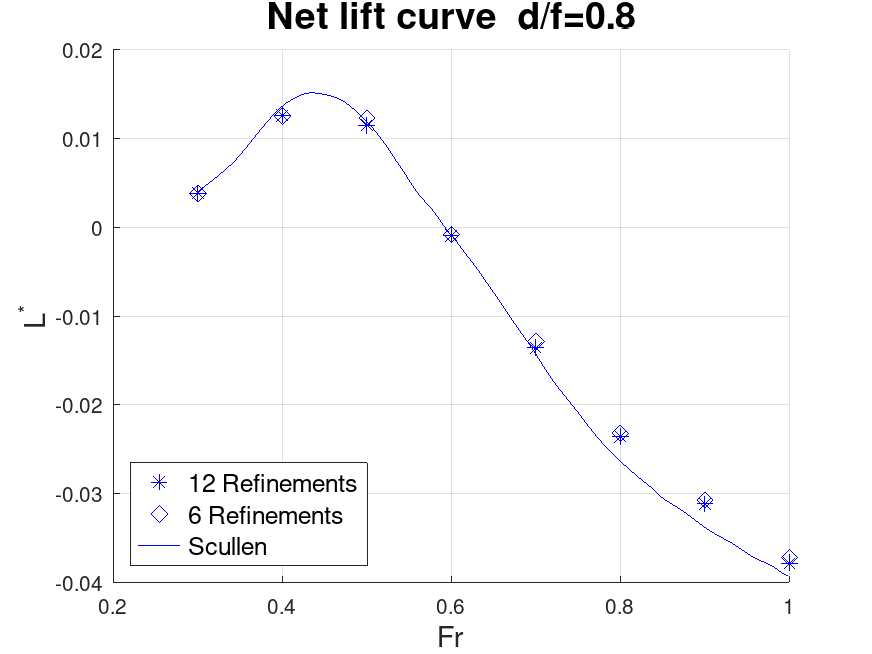}
\includegraphics[width=0.49\linewidth]{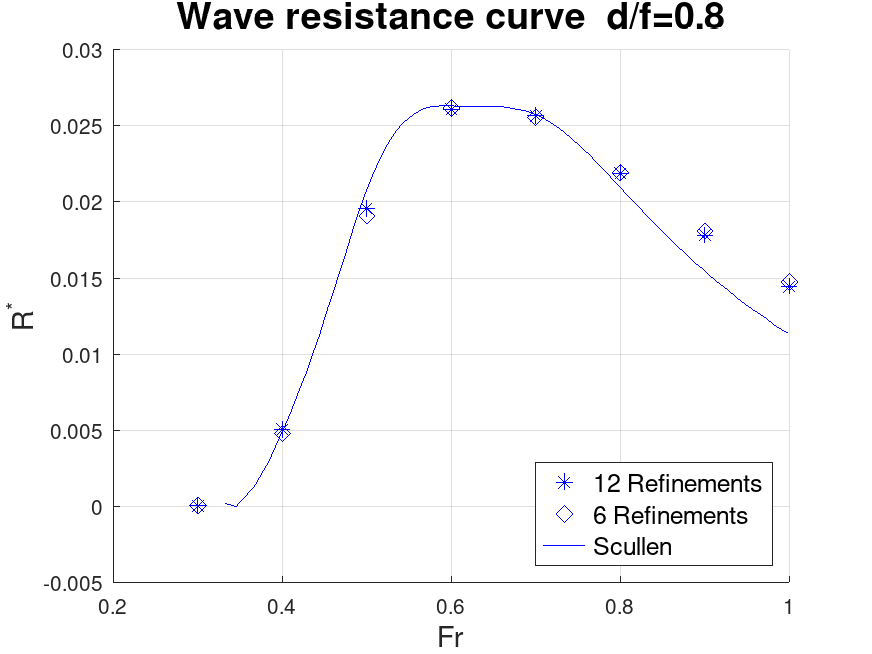}
\caption{Fully immersed ellipsoid non dimensional net lift $L^*=(L-L_0)/L_0$ (left) and resistance
         $R^*=(R)/L_0$ (right) as a function of Froude
         number ${Fr}$. For reference, the continuous lines refers to results obtained by
         Scullen for ellipsoid depth to diameter ratio $d/f = 0.8$. The results obtained at the 6th
         and 12th adaptive refinement iterations are indicated by diamonds and asterisks,
         respectively.\label{fig:ForceCompareAcrossRefs}}
\end{figure}
\end{center}

To provide the reader with an extensive assessment the accuracy of the model proposed, we
compare the nondimensional forces computed with the present method using 7 refinement cycles
against similar results reported in ~\cite{Scullen}. In such work, Scullen made use of a steady potential
flow solver with fully nonlinear free surface boundary conditions which represent a
combination of the null pressure and null pressure total derivative requirement. 
Figure~\ref{fig:scullenCompareLift} compares the non dimensional net lift values obtained in this
work (indicated by asterisks) against the corresponding results obtained by Scullen (solid lines).
The different curves in the plot represent non dimensional net lift as a function of Fr obtained
imposing different values of the $\frac{d}{f}$ ratio. The results show good agreement with Scullen
data throughout the Fr and $\frac{d}{f}$ ranges tested. The plot also indicates that the most
appreciable differences are observed for higher Fr and $\frac{d}{f}$ values, where the present
method lift is consistently higher than its reference counterpart. This might be a result
of the different --- and possibly less dissipative --- BEM formulation in which Rankine sources
coincide with collocation nodes, and special singular quadrature is used. Such formulation
might result in higher free surface nodes displacements which bring the free surface closer to the
spheroid in the rear region, leading to increased local speeds associated with upwards
suction. 

\begin{center}
\begin{figure}
\includegraphics[width=0.9\linewidth]{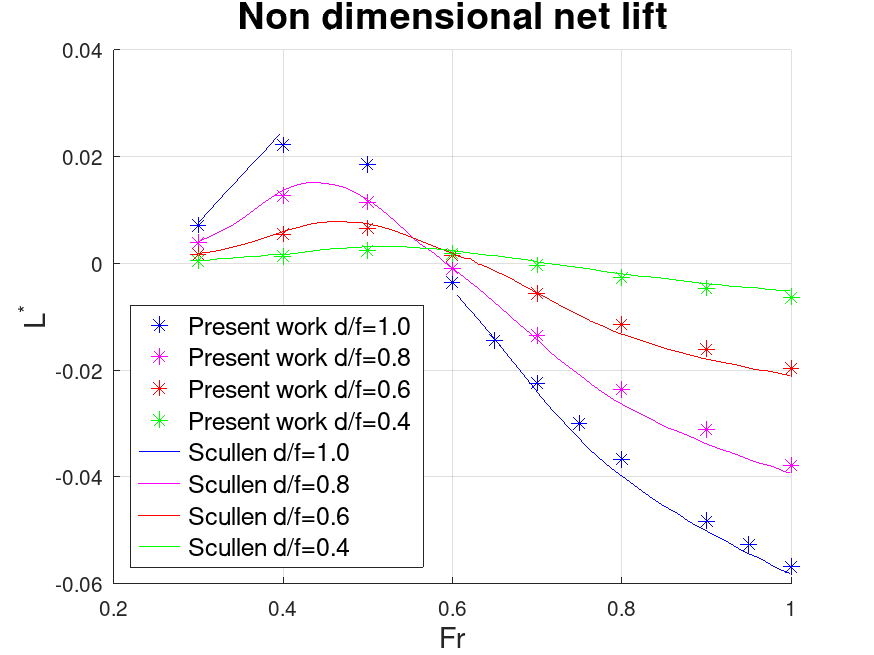}
\caption{Fully immersed ellipsoid non dimensional net lift $L^*=(L-L_0)/L_0$ as a function of
         Froude number Fr. The differently colored continuous lines refer to results obtained
         by Scullen for values of ellipsoid depth to
         diameter ratio $d/f = \textcolor{blue}{1.0}, \textcolor{magenta}{0.8}, \textcolor{red}{0.6},
         \textcolor{green}{0.4}$. The results obtained in this work are indicated by asterisks with
         corresponding colors.\label{fig:scullenCompareLift}}
\end{figure}
\end{center}

Figure~\ref{fig:scullenCompareResistance} displays a comparison of non dimensional wave resistance
values obtained in this work (indicated by asterisks) against the corresponding results obtained
by Scullen (solid lines). The different curves in the plot represent non dimensional resistance
as a function of Fr obtained imposing different values of the $\frac{d}{f}$ ratio. Also in this
case, the results seem in good agreement with their reference literature counterparts across
all the range of Fr and $\frac{d}{f}$ values considered. Once again, the most relevant differences
are observed in correspondence with the highest Fr and $\frac{d}{f}$ tested, where the wave drag
predicted by the present method is consistently and sensibly higher than the value reported by
Scullen. This can once again be ascribed to the less dissipative nature of the solver proposed,
which predicts higher surface displacements bringing the free surface closer to the hull in the
stern region. This generates higher suction and, in turn, an additional drag due to lower
pressures in the stern region.

\begin{center}
\begin{figure}
\includegraphics[width=0.9\linewidth]{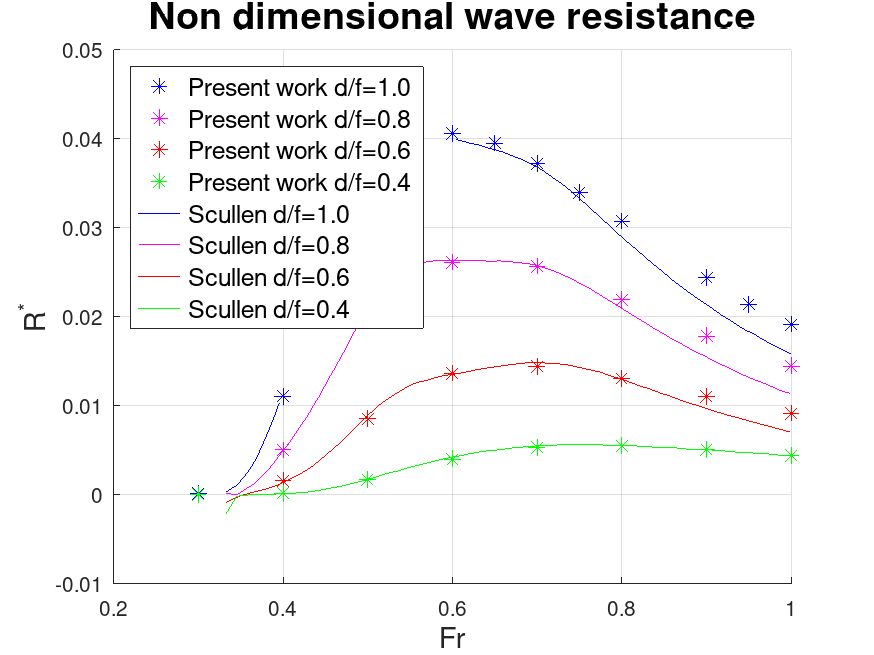}
\caption{Fully immersed ellipsoid non dimensional wave resistance $R^*=(R)/L_0$ as a function of
         Froude number Fr. The differently colored continuous lines refer to results obtained
         by Scullen for values of ellipsoid depth to diameter
         ratio $d/f = \textcolor{blue}{1.0}, \textcolor{magenta}{0.8}, \textcolor{red}{0.6},
         \textcolor{green}{0.4}$. The results obtained in this work are indicated by asterisks with
         corresponding colors. \label{fig:scullenCompareResistance}}
\end{figure}
\end{center}

We must finally remark that as Figures \ref{fig:scullenCompareLift} and \ref{fig:scullenCompareResistance}
suggest, some solutions have been obtained with the present method in regions where previous methods
seem not reach convergence. This is particularly evident for $\frac{d}{f}=1$ and Fr$=[0.4,0.6]$,
for which no solution is reported by Scullen. This should not surprise, as in such conditions the
deepest wave trough is located slightly above the ellipsoid stern region, and almost contact
is reached between the free surface and the hull. In the framework of the current method
implementation, no contact between free surface and disturbing body is considered. Moreover,
the non desingularized method used has currently no mechanism in place to apply singular
quadrature on the hull stern surface and account for the very close Rankine sources in the free
surface trough, and vice versa. For such reason, also in this work some simulations failed to
reach convergence at the last refinement levels. Thus, we must report that the ability to predict
hydrodynamic forces in such condition is a result of the adaptive refinement approach adopted
rather than a product of a superior robustness. Yet, being able to obtain a solution with a locally
less refined grid, gives us the opportunity to obtain a viable drag and lift estimate also in such
difficult test case.

\subsubsection{Results of unsteady case with steady regime solution}

The numerical results of the unsteady test cases described in Section \ref{sec:fake_unsteady_experiments}
are now taken into consideration. Figure \ref{fig:steadyVsunsteadyDrag} presents a plot of the
underwater ellipsoid wave resistance force absolute value as a function of time. The blue, green and
magenta continuous lines refer to the test cases in which the flow asymptotic velocity $\Ub_\infty$
reaches the target value after sinusoidal ramps of $0.75\,s$, $7.5\,s$ and $15\,s$, respectively. The
diagram clearly shows that all the test cases considered gradually reach the same constant regime
solution. For the purposes of the present work, it is important to point out that the
wave resistance value associated with such a common regime solution is identical to the one obtained with
the steady state simulation on the same computational grid, denoted by the dashed red line in the plot.
On one hand, this should not be a surprise, as the nonlinear problems solved for the steady and unsteady
solver are the same if $t\rightarrow\infty$, as illustrated in Section \ref{sec:steady}. On the
other hand, it must be stressed that obtaining a potential flow solver with fully nonlinear free surface
treatment and the ability to compute both transient and steady solutions, is one of the main objectives
of this work. A further look at the time evolution plot for the three test cases, shows that, as expected,
the test case with faster dynamics (blue line) results in a higher peak resistance associated with the
added mass contribution induced by the increased initial acceleration. As a consequence of this, also the
rebound resistance local minimum following the initial water acceleration past the hull, is more intense
in the $0.75\,s$ test case.

\begin{center}
\begin{figure}
\includegraphics[width=0.9\linewidth]{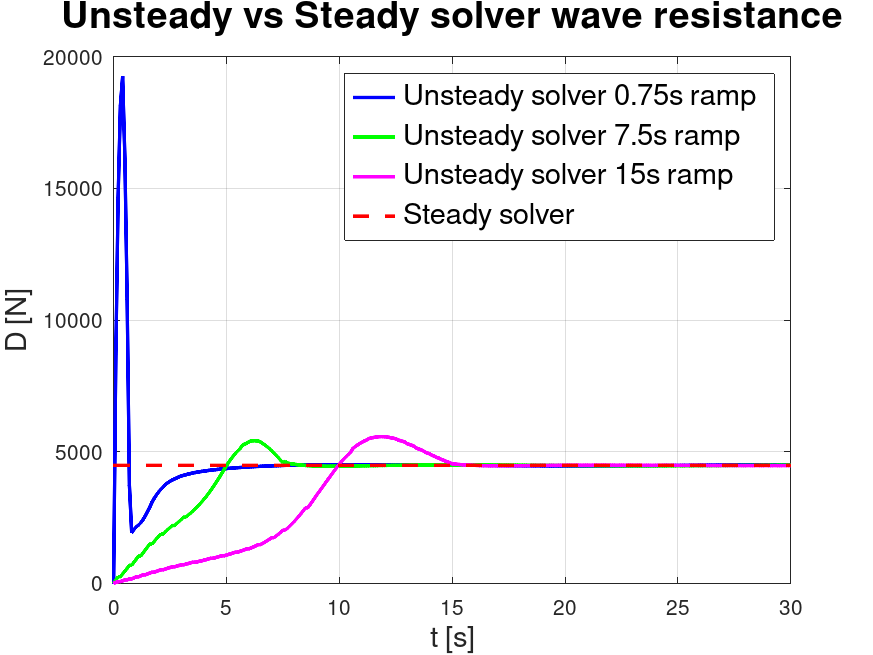}
\caption{Wave resistance as a function of time for the fully immersed ellipsoid
         with $d/f = 0.8$ advancing at a regime speed corresponding to $\text{Fr}=0.8$.
         The red dashed curve represents the result obtained with the steady solver
         after 3 adaptive refinement cycles. The unsteady solver simulation results
         obtained with a $0.75\,s$, $7.5\,s$ and $15\,s$ are displayed by the
         blue, green and magenta curve, respectively.\label{fig:steadyVsunsteadyDrag}}
\end{figure}
\end{center}

Figure \ref{fig:steadyVsunsteadyLifts} depicts absolute values of hydrodynamic lift acting on the
immersed ellipsoid, as a function of time. In the diagram, the lines color are associated
to the $0.75\,s$, $7.5\,s$ and $15\,s$ ramp test cases, in the same way reported for the
wave resistance plot. The lift plots substantially confirm what previously observed analyzing
the resistance results. Also in this case, all the unsteady flow solver results appear to
converge to a common steady state lift value, which coincides with the value resulting from
the steady solver simulation on the same grid.

\begin{center}
\begin{figure}
\includegraphics[width=0.9\linewidth]{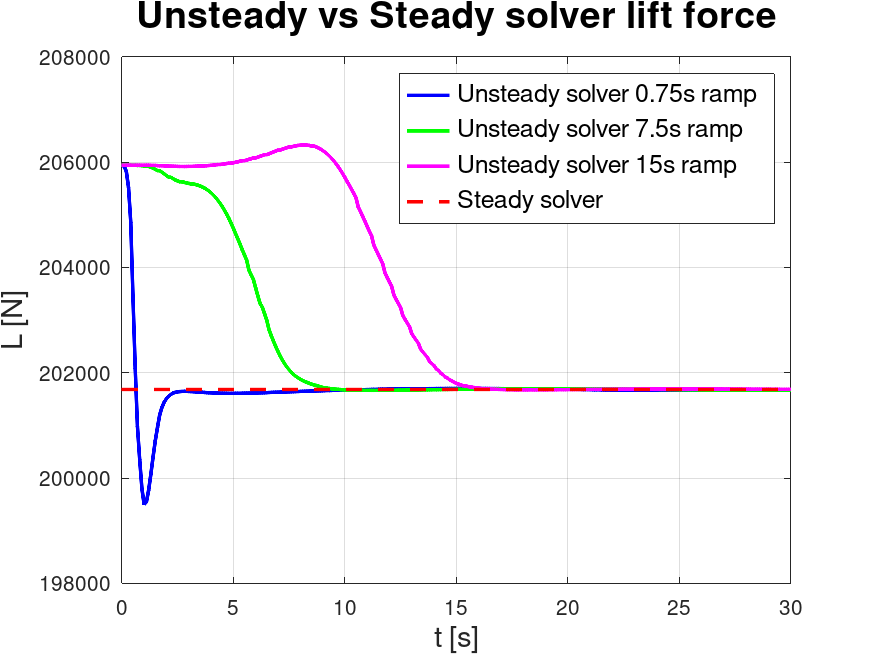}
\caption{Hydrodynamic lift as a function of time for the fully immersed ellipsoid
         with $d/f = 0.8$ advancing at a regime speed corresponding to $\text{Fr}=0.8$.
         The red dashed curve represents the result obtained with the steady solver
         after 3 adaptive refinement cycles. The unsteady solver simulation results
         obtained with a $0.75\,s$, $7.5\,s$ and $15\,s$ are displayed by the
         blue, green and magenta curve, respectively.\label{fig:steadyVsunsteadyLifts}}
\end{figure}
\end{center}

\subsubsection{Results of unsteady case with unsteady regime solution}

Figure \ref{fig:fully_unsteady_frames} illustrates two different time steps of the
numerical solution obtained for the unsteady test cases described in Section \ref{sec:true_unsteady_experiments}. The
unsteady flow simulation of $100\,$s, corresponding to approximately 20 incident wave periods, lasted 5 to 6 hours.
The plots clearly show nonlinear interaction between the unsteady incident waves arriving towards the
ellipsoid, and the Kelvin wake pattern, which is located in a steady position in the moving reference frame.  
The two images in the figure show instants in which a thorough (left plot) and a crest (right plot)
of the incident wave field pass on top of the spheroid. Both images also include a view 
of the numerical domain up until its end. The flow approaching the truncation boundary of the numerical domain
appears to converge to the asymptotic Airy wave imposed, suggesting that the numerical damping zone implemented
successfully absorbs waves potentially reflected in the flow field.

\begin{center}
\begin{figure}
\begin{tabular}{c c}
\includegraphics[width=0.475\linewidth]{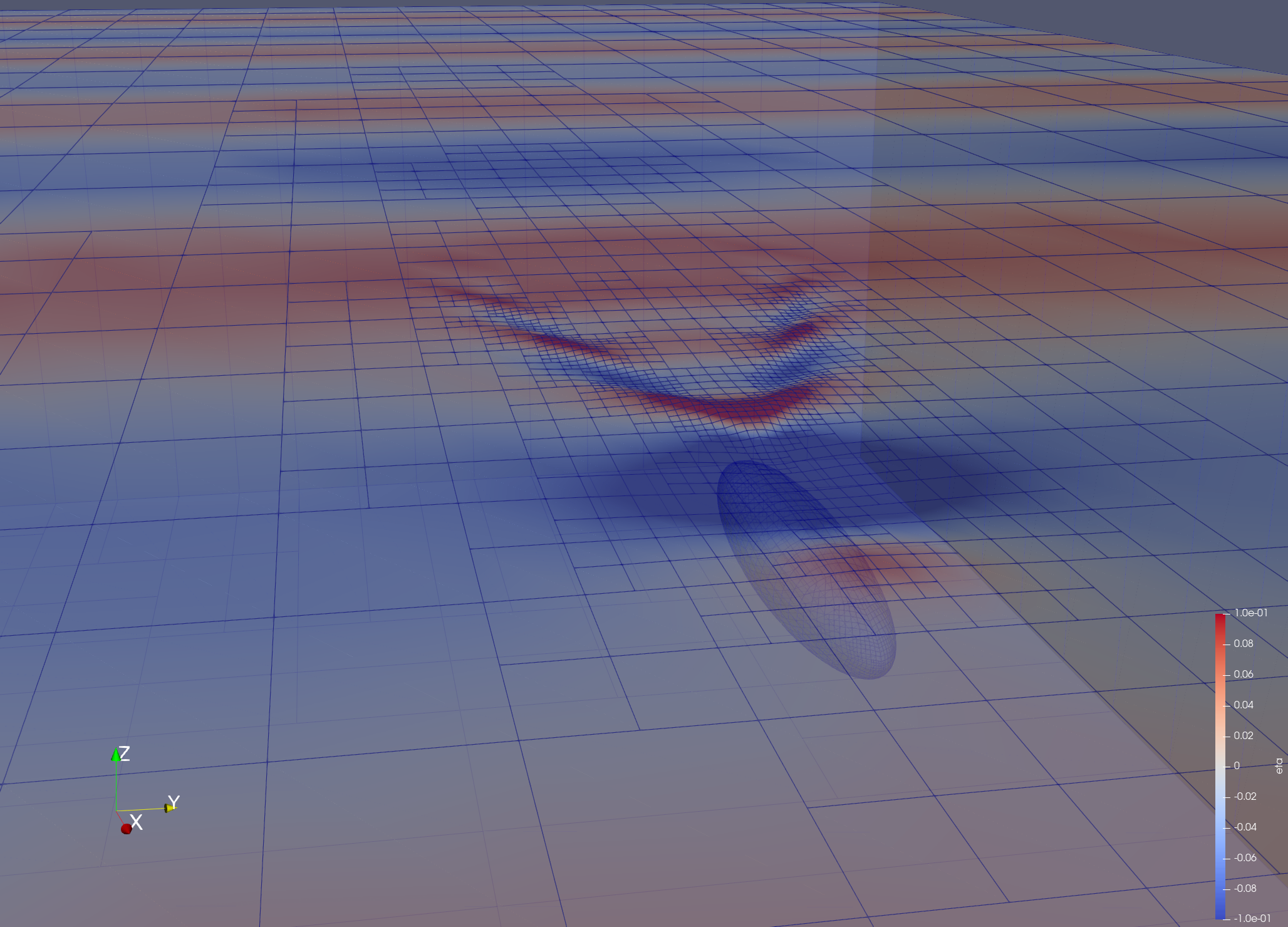}
&
\includegraphics[width=0.475\linewidth]{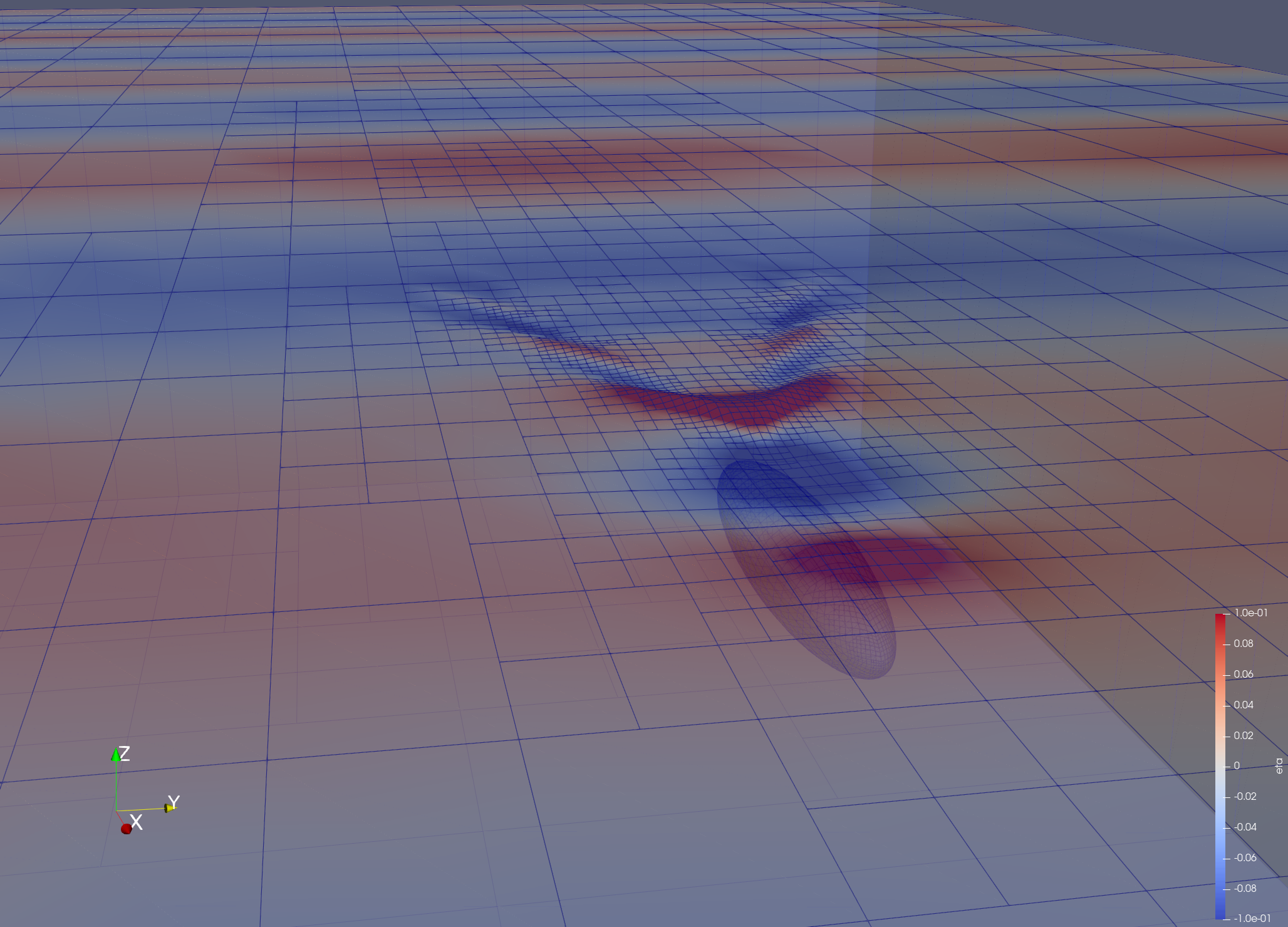}
\end{tabular}
\caption{Two different time frames depicting the water free surface surrounding the fully
         submerged ellipsoid advancing at $\text{Fr}=0.3$ in head waves with
         asymptotic wavelength $\lambda=4L$ and amplitude $a=0.02\,$m. The free surface
         is colored based on contours of the water elevation field. Note that the water elevation
         peak in the flow field is approximately $0.21\,$m. However, in the plot
         the water elevation contour scale range is limited to $[-0.1\, \text{m},0.1\, \text{m}]$
         so as to better show features of the incident wave pattern.
         A full video of the present simulation is available to the interested reader
         at the URL address \protect\url{https://youtu.be/t4XSTJ4RxFU}.
         \label{fig:fully_unsteady_frames}}
\end{figure}
\end{center}

To provide a visual evaluation of the incident wavelength observed in the nonlinear basin,
Figure \ref{fig:fully_unsteady_top} shows a vertical view of the water elevation field
around the hull. In particular, the diagram shows that the distance between the
two crests approaching the ellipsoid is a good approximation of the asymptotic
value $\lambda=4L$ imposed. The wave amplitude in instead higher than the value imposed in the
asymptotic potential, and settles for values between $0.04\,$m and $0.05\,$m. This could be 
obviously due to numerical error, or to blocking effects due to the presence of the ellipsoid
in the water channel. A further possible explanation is that the discrepancy is due to a
natural shape difference between the linear waves imposed at a distance from the hull and
the nonlinear ones developing once they enter in the domain. Further investigations will
be carried out to obtain a better control of the water elevation amplitude.

\begin{center}
\begin{figure}
\begin{center}
\includegraphics[width=0.75\linewidth]{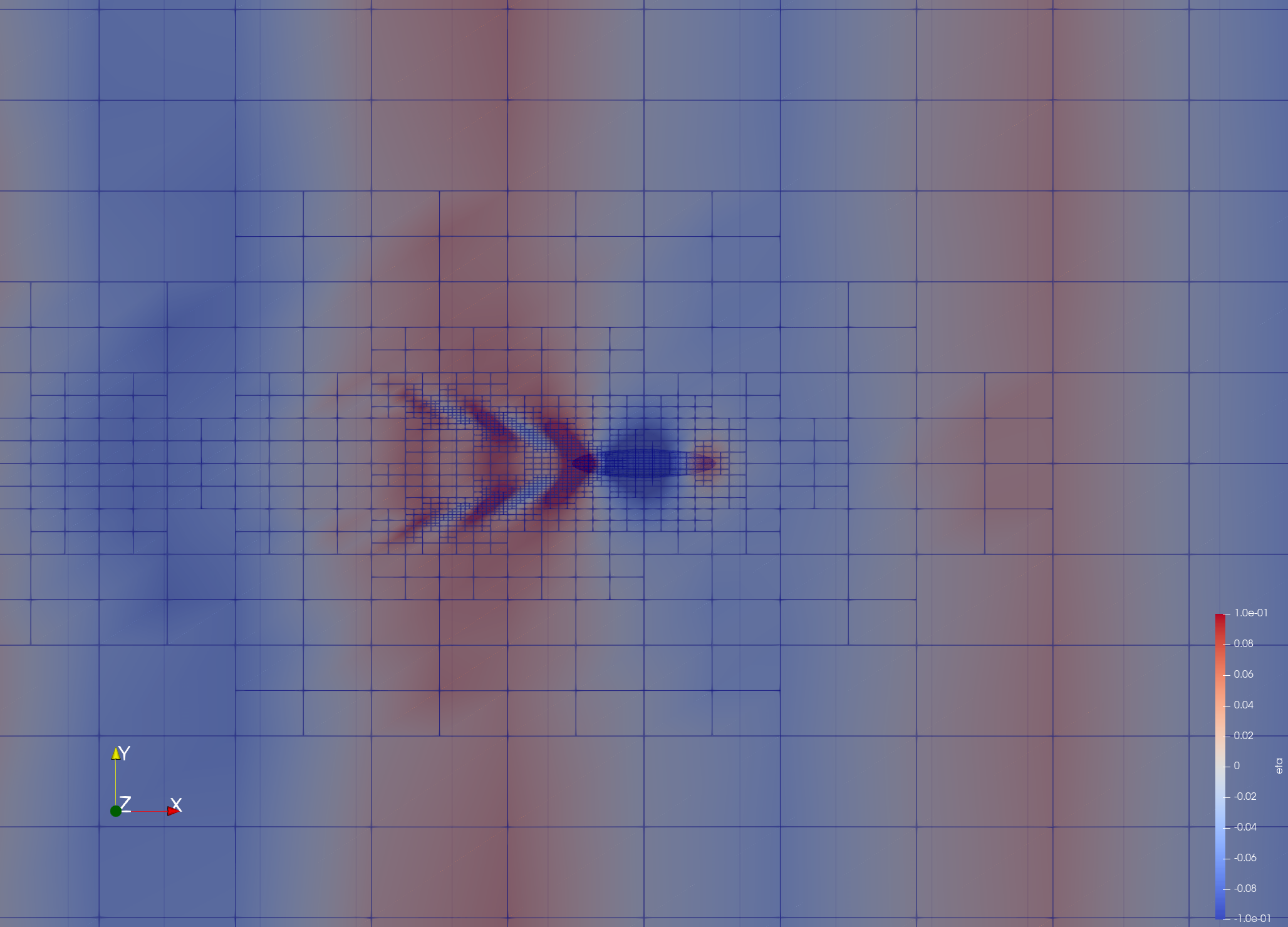}
\caption{A vertical view of the wave field surrounding the fully
         submerged ellipsoid advancing at $\text{Fr}=0.3$ in head waves with
         asymptotic wavelength $\lambda=4L$ and amplitude $a=0.02\,$m. The free surface
         is colored based on contours of the water elevation field.
         Note that the water elevation
         peak in the flow field is approximately $0.21\,$m. However, in the plot
         the water elevation contour scale range is limited to $[-0.1\, \text{m},0.1\, \text{m}]$
         so as to better show features of the incident wave pattern.
         \label{fig:fully_unsteady_top}}
\end{center}
\end{figure}
\end{center}

\begin{center}
\begin{figure}[t]
\includegraphics[width=0.9\linewidth]{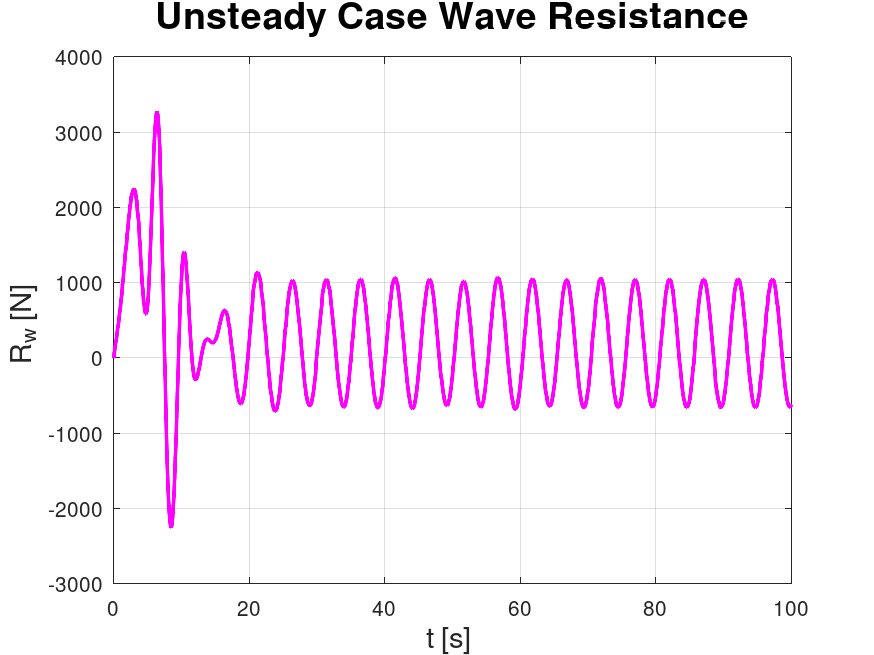}
\caption{Wave resistance as a function of time for the fully immersed ellipsoid
         with $d/f = 0.75$ advancing at a regime speed corresponding to $\text{Fr}=0.3$
         in an Airy wave field with wavelength $\lambda=40\,$m and $a=0.02\,$m.
         \label{fig:fully_unsteadyDrag}}
\end{figure}
\end{center}

Figure \ref{fig:fully_unsteadyDrag} shows the time history of the ellipsoid wave resistance
measured during the unsteady flow simulation. The plot shows that the simulation 
appears stable for its entire duration of $100\,$s. The unsteady flow settles for
a periodic regime solution, suggesting once again that the damping strategy used is
able to successfully prevent wave energy from being reflected back into the domain, which
leads to stable simulations. The wave resistance period as observed in the plot is slightly
higher than $5\,$s, which is close to the expected value $T = 5.0616\, s$ imposed in the
asymptotic flow. By a qualitative perspective, the time history plot suggests that wave
resistance is dominated by its added mass contribution, as the resistance values follow the
local wave field acceleration, and even become negative when the
ellipsoid encounters waves accelerating towards its forward motion direction.  The average wave
resistance value computed from the last 6 full periods
of the simulation (lasting approximately $30\,$s) is $\overline{R}=181.66\,$N. When compared
to the steady state wave resistance $R_{SS}=-178.23\,$N, this results in a positive added resistance
value due to waves $\Delta R_w=3.43\, $N, corresponding to approximately 1.9\% of the steady resistance
value. A similar test case has been proposed by Papanikolaou and Liu in \cite{PapanLiu2010}, which report
(see Figure 9 in such paper) non dimensional added resistance $\Delta R_w^{REF}/(\rho g a^2 d^2/L)=0.145$ computed by means of a
linearized free surface model featuring a pulsating and traveling { Green's function} taking into account
forward speed effects. The corresponding dimensional value is $\Delta R_w^{REF} = 0.23\,$N if
a wave amplitude $a=0.02\ $m is considered. Alternatively, if the wave amplitude considered is the
one effectively observed in the basin ($a=0.05\,$m), the dimensional value based on Papanikolaou and Liu
is $\Delta R_w^{REF} = 1.42\,$N. Both values are considerably lower than the one computed in the present
simulation campaign. The difference can be in principle be associated to the different models used,
as it is by all means possible that a nonlinear free surface model can result in higher added wave resistance.
However, further investigations will be carried out in the near future to obtain a full characterization
of the accuracy of added wave resistance predicted by the present model.

\section{Conclusions and future perspectives}

This work presented a novel formulation of the fully nonlinear free surface boundary conditions
which complement the Laplace equation in numerical towing tank based on unsteady potential flow theory.
The main advantage of the unsteady free surface model discussed, is that it can lead to steady
state solutions once the --- null --- time derivatives are eliminated by the discretized problem.
Such a feature, as discussed, is not common in fully nonlinear potential flow solvers available
in the literature and in the commercial software market. Numerical results presented
confirm that the steady and unsteady solvers result in the same solution for long time integration,
and that the steady solver solutions are in close comparison with classical steady nonlinear free surface
potential solvers.

A possible immediate future work perspective, as mentioned in the text, is carrying out a simulation
campaign to fully characterize the quality of the waves introduced in the model making use of the asymptotic potential.
In particular, investigating whether different incident waves directions and multichromatic waves can be
successfully considered should result in interesting perspective works. In addition, other work should focus on
implementing the CAD interface \cite{Mola2014} and Fluid Structure interaction solver for rigid
ships \cite{molaIsope2016} also in the present software. This would result in
a versatile and effective numerical towing tank for ship hydrodynamics problems. Also adding the effect
of lifting surfaces as in \cite{saccoMaster2017} would result in including the effects of
hull appendages in the model.

Finally, we must point out that the free surface boundary condition in ALE form used in this work
does not depend on the assumption that the free surface $\eta$ is a single valued Cartesian function.
In fact, it is only through Equation \eqref{eq:v_grid} that such constraint is introduced in the
system, whereas Equations \eqref{eq:kin_fs_condition_non_pen} and \eqref{eq:dyn_fs_condition_semil}
are the ALE version of the non penetration and dynamic free surface condition, and can in principle
work with any grid velocity field $\vb$. Thus, future work will investigate the possibility of reproducing
steep or overturning waves in the present formulation, through a wiser $\vb$ choice.

\section{Acknowledgments}

We acknowledge the support by the European Commission H2020 ARIA (Accurate ROMs for Industrial Applications, GA 872442) project, by MIUR (Italian Ministry for Education University and Research) and PRIN "Numerical Analysis for Full and Reduced Order Methods for Partial Differential Equations" (NA-FROMPDEs) project, by the European Research Council Consolidator Grant Advanced Reduced Order Methods with Applications in Computational Fluid Dynamics-GA 681447, H2020-ERC COG 30 2015 AROMA-CFD, by FSE Galicia 2014-2020 and Xunta de Galicia under grant ED481A-2018/212, by FEDER, Ministerio de Economía, Industria y Competitividad-AEI research project MTM2017-86459-R, Ministerio de Ciencia e Innovación through the research project  PID2021-122625OB-I00 and by Xunta de Galicia (Spain) research project GI-1563 ED431C 2021/15.

\newpage
\appendix 
\addcontentsline{toc}{section}{Appendices}

\section{Summary of numerical scheme features}

\begin{center}
\begin{table}[h!]
    \hspace*{3.5cm}
    \begin{sideways}
    \tiny
    \begin{tabular}{c|c|c|c|c} 
          \textbf{Type}   & \textbf{Type of}      & \textbf{Treats time}        & \textbf{Able to reach steady} & \textbf{Leads to } \\
          \textbf{of}     & \textbf{time}         & \textbf{dependence and}     & \textbf{state in hull }       & \textbf{numerically stable} \\
          \textbf{solver} & \textbf{derivatives}  & \textbf{moving domains}     & \textbf{frame of reference}   & \textbf{steady problem} \\
      \hline
                                 &                     &     &     &     \\[4mm]
      
      Lagrangian MEL \cite{grilli2001,longuet-higginsCokeletMEL76}            & Lagrangian          & {\color{green}{Yes}} & {\color{red}{No}}  & {\color{red}{/}}   \\[4mm]
      
      semi-Lagrangian MEL \cite{beck1994,scorpioPhD1997}       & ALE                 & {\color{green}{Yes}} & {\color{red}{No}}  & {\color{red}{/}}   \\[4mm]
      
      stabilized semi-Lagrangian \cite{waveBem} & ALE                 & {\color{green}{Yes}} & {\color{green}{Yes}} & {\color{red}{No}}  \\[4mm]
      
      steady state solvers \cite{ravenPhD1998,JansonPhD1997,Scullen}      & Eulerian (set to 0) & {\color{red}{No}}    & {\color{green}{Yes}} & {\color{green}{Yes}} \\[4mm]
      
      current model              & ALE                 & {\color{green}{Yes}} & {\color{green}{Yes}} & {\color{green}{Yes}} \\[4mm]
    \end{tabular}
    \end{sideways}
    \caption{Comparison between the most common numerical approaches to the nonlinear free surface problem. \label{tab:techniques_comparison}}
\end{table}

\end{center}
\normalsize

\FloatBarrier

\vspace{1cm}
\section{Resolution algorithm flow chart}

\begin{center}
\begin{figure}[h!]
\centerline{
  \ifpdf
  \resizebox{0.75\textwidth}{!}{
    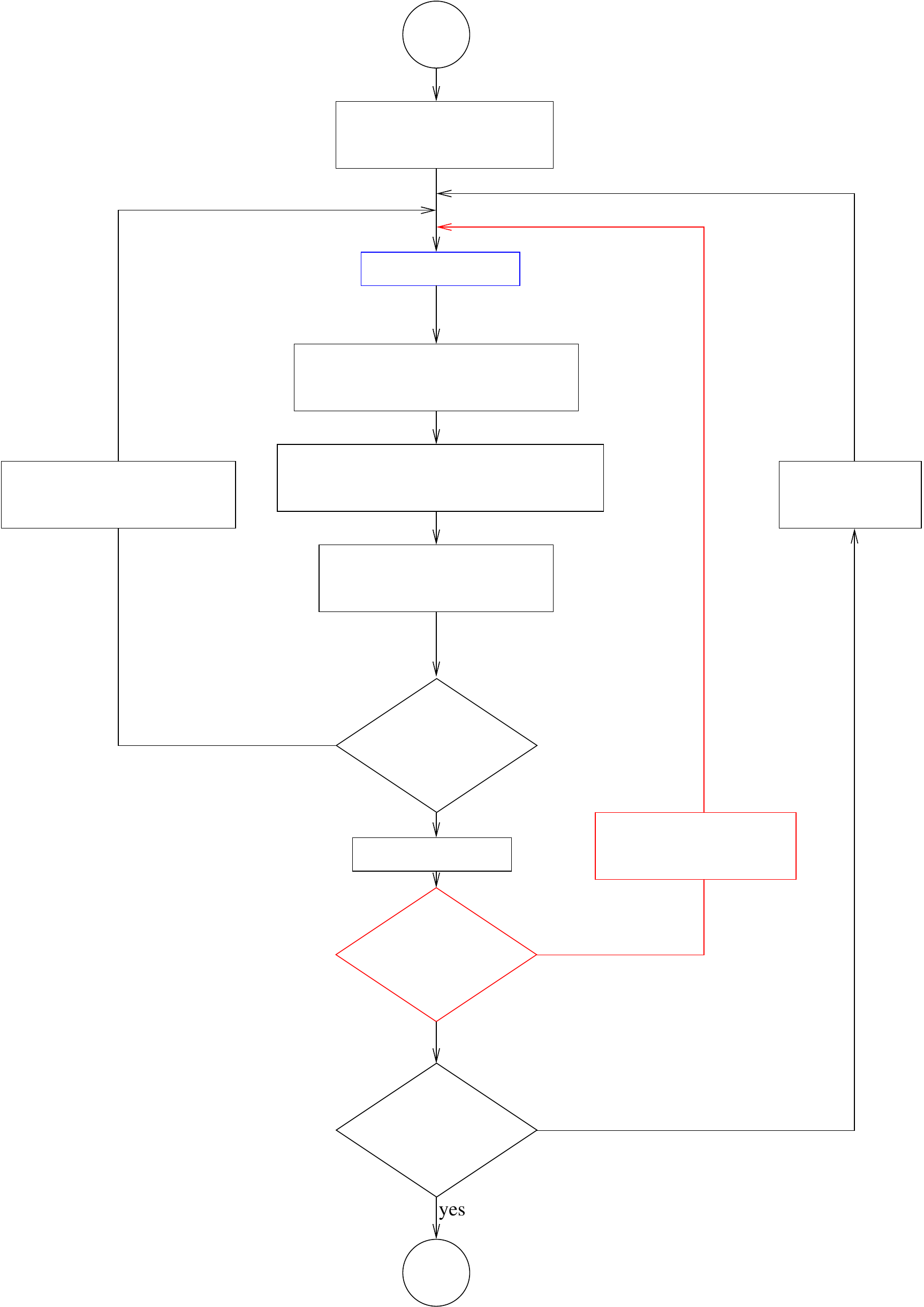
  }
  \else
  \resizebox{0.75\textwidth}{!}{
    \input{./figures/flow_chart.pstex_t}
  }
  \fi
}
\caption{A flow chart representing the numerical resolution algorithms. Note that the
         red portions of the chart are only executed in case the unsteady solver is used, while the
         blue parts are only executed for steady runs.\label{fig:flow_chart}}
\end{figure}
\end{center}

\FloatBarrier

\bibliography{mybibfile}

\end{document}

%% file: figures/dominio_barca.pdftex_t
\begin{picture}(0,0)%
\includegraphics{figures/dominio_barca.pdf}%
\end{picture}%
\setlength{\unitlength}{2693sp}%
\begingroup\makeatletter\ifx\SetFigFont\undefined%
\gdef\SetFigFont#1#2#3#4#5{%
  \reset@font\fontsize{#1}{#2pt}%
  \fontfamily{#3}\fontseries{#4}\fontshape{#5}%
  \selectfont}%
\fi\endgroup%
\begin{picture}(8858,7023)(2303,-7842)
\put(5639,-3134){\makebox(0,0)[lb]{\smash{{\SetFigFont{8}{9.6}{\familydefault}{\mddefault}{\updefault}{\color[rgb]{0,0,0}$y$}%
}}}}
\put(7469,-3314){\makebox(0,0)[lb]{\smash{{\SetFigFont{8}{9.6}{\familydefault}{\mddefault}{\updefault}{\color[rgb]{0,0,0}$x$}%
}}}}
\put(6569,-2789){\makebox(0,0)[lb]{\smash{{\SetFigFont{8}{9.6}{\familydefault}{\mddefault}{\updefault}{\color[rgb]{0,0,0}$z$}%
}}}}
\put(2543,-1871){\makebox(0,0)[lb]{\smash{{\SetFigFont{8}{9.6}{\familydefault}{\mddefault}{\updefault}{\color[rgb]{0,0,0}$\Gamma^{fs}$}%
}}}}
\put(9247,-7037){\makebox(0,0)[lb]{\smash{{\SetFigFont{8}{9.6}{\familydefault}{\mddefault}{\updefault}{\color[rgb]{0,0,0}$\Gamma^{\infty}$}%
}}}}
\put(4561,-1501){\makebox(0,0)[lb]{\smash{{\SetFigFont{8}{9.6}{\familydefault}{\mddefault}{\updefault}{\color[rgb]{0,0,0}$\Gamma^h$}%
}}}}
\put(6661,-7351){\makebox(0,0)[lb]{\smash{{\SetFigFont{8}{9.6}{\familydefault}{\mddefault}{\updefault}{\color[rgb]{0,0,0}$\Gamma^{b}$}%
}}}}
\put(3106,-7216){\makebox(0,0)[lb]{\smash{{\SetFigFont{8}{9.6}{\familydefault}{\mddefault}{\updefault}{\color[rgb]{0,0,0}$\Gamma^{i}$}%
}}}}
\put(2926,-5686){\makebox(0,0)[lb]{\smash{{\SetFigFont{8}{9.6}{\familydefault}{\mddefault}{\updefault}{\color[rgb]{0,0,0}$\nablab\phi_{\infty}$}%
}}}}
\end{picture}%

%% file: figures/spheroid_domain.pdftex_t
\begin{picture}(0,0)%
\includegraphics{figures/spheroid_domain.pdf}%
\end{picture}%
\setlength{\unitlength}{4144sp}%
\begingroup\makeatletter\ifx\SetFigFont\undefined%
\gdef\SetFigFont#1#2#3#4#5{%
  \reset@font\fontsize{#1}{#2pt}%
  \fontfamily{#3}\fontseries{#4}\fontshape{#5}%
  \selectfont}%
\fi\endgroup%
\begin{picture}(9899,3904)(699,-7103)
\put(5041,-3391){\makebox(0,0)[lb]{\smash{{\SetFigFont{10}{12.0}{\familydefault}{\mddefault}{\updefault}{\color[rgb]{0,0,0}$z$}%
}}}}
\put(5806,-4111){\makebox(0,0)[lb]{\smash{{\SetFigFont{10}{12.0}{\familydefault}{\mddefault}{\updefault}{\color[rgb]{0,0,0}$x$}%
}}}}
\put(1756,-5641){\makebox(0,0)[lb]{\smash{{\SetFigFont{14}{16.8}{\familydefault}{\mddefault}{\updefault}{\color[rgb]{0,0,0}$15L$}%
}}}}
\put(7201,-4471){\makebox(0,0)[lb]{\smash{{\SetFigFont{14}{16.8}{\familydefault}{\mddefault}{\updefault}{\color[rgb]{0,0,0}$f$}%
}}}}
\put(8596,-5641){\makebox(0,0)[lb]{\smash{{\SetFigFont{14}{16.8}{\familydefault}{\mddefault}{\updefault}{\color[rgb]{0,0,0}$15L$}%
}}}}
\put(5176,-5641){\makebox(0,0)[lb]{\smash{{\SetFigFont{14}{16.8}{\familydefault}{\mddefault}{\updefault}{\color[rgb]{0,0,0}$L$}%
}}}}
\put(2386,-4921){\makebox(0,0)[lb]{\smash{{\SetFigFont{14}{16.8}{\familydefault}{\mddefault}{\updefault}{\color[rgb]{0,0,0}$d = \displaystyle{\frac{2L}{5}}$}%
}}}}
\put(1036,-4426){\makebox(0,0)[lb]{\smash{{\SetFigFont{14}{16.8}{\familydefault}{\mddefault}{\updefault}{\color[rgb]{0,0,0}$\nablab\phi_\infty(t)$}%
}}}}
\put(9811,-6586){\makebox(0,0)[lb]{\smash{{\SetFigFont{14}{16.8}{\familydefault}{\mddefault}{\updefault}{\color[rgb]{0,0,0}$5L$}%
}}}}
\end{picture}%

%% file: figures/flow_chart.pdftex_t
\begin{picture}(0,0)%
\includegraphics{figures/flow_chart.pdf}%
\end{picture}%
\setlength{\unitlength}{3947sp}%
\begingroup\makeatletter\ifx\SetFigFont\undefined%
\gdef\SetFigFont#1#2#3#4#5{%
  \reset@font\fontsize{#1}{#2pt}%
  \fontfamily{#3}\fontseries{#4}\fontshape{#5}%
  \selectfont}%
\fi\endgroup%
\begin{picture}(8274,11714)(3139,-9143)
\put(8551,-5191){\makebox(0,0)[lb]{\smash{{\SetFigFont{12}{14.4}{\familydefault}{\mddefault}{\updefault}{\color[rgb]{1,0,0}and update time}%
}}}}
\put(10351,-2011){\makebox(0,0)[lb]{\smash{{\SetFigFont{12}{14.4}{\familydefault}{\mddefault}{\updefault}{\color[rgb]{0,0,0}refine grid}%
}}}}
\put(6470,-5133){\makebox(0,0)[lb]{\smash{{\SetFigFont{10}{12.0}{\familydefault}{\mddefault}{\updefault}{\color[rgb]{0,0,0}Output solution}%
}}}}
\put(6826,-5826){\makebox(0,0)[lb]{\smash{{\SetFigFont{12}{14.4}{\familydefault}{\mddefault}{\updefault}{\color[rgb]{1,0,0}Time}%
}}}}
\put(6989,-6003){\makebox(0,0)[lb]{\smash{{\SetFigFont{12}{14.4}{\familydefault}{\mddefault}{\updefault}{\color[rgb]{1,0,0}=}%
}}}}
\put(6631,-6185){\makebox(0,0)[lb]{\smash{{\SetFigFont{12}{14.4}{\familydefault}{\mddefault}{\updefault}{\color[rgb]{1,0,0}Max time?}%
}}}}
\put(6553,-7799){\makebox(0,0)[lb]{\smash{{\SetFigFont{12}{14.4}{\familydefault}{\mddefault}{\updefault}{\color[rgb]{0,0,0}Max cycles?}%
}}}}
\put(8287,-7509){\makebox(0,0)[lb]{\smash{{\SetFigFont{12}{14.4}{\familydefault}{\mddefault}{\updefault}{\color[rgb]{0,0,0}no}%
}}}}
\put(6457,121){\makebox(0,0)[lb]{\smash{{\SetFigFont{10}{12.0}{\familydefault}{\mddefault}{\updefault}{\color[rgb]{0,0,1}Set $\dot{\yb}=0, t=\infty$}%
}}}}
\put(6226,1439){\makebox(0,0)[lb]{\smash{{\SetFigFont{12}{14.4}{\familydefault}{\mddefault}{\updefault}{\color[rgb]{0,0,0}Prepare initial grid}%
}}}}
\put(6226,1184){\makebox(0,0)[lb]{\smash{{\SetFigFont{12}{14.4}{\familydefault}{\mddefault}{\updefault}{\color[rgb]{0,0,0}and initialize solution}%
}}}}
\put(5851,-736){\makebox(0,0)[lb]{\smash{{\SetFigFont{12}{14.4}{\familydefault}{\mddefault}{\updefault}{\color[rgb]{0,0,0}Assemble BEM system (42)}%
}}}}
\put(6151,-961){\makebox(0,0)[lb]{\smash{{\SetFigFont{12}{14.4}{\familydefault}{\mddefault}{\updefault}{\color[rgb]{0,0,0}and store Derivatives}%
}}}}
\put(6076,-1861){\makebox(0,0)[lb]{\smash{{\SetFigFont{12}{14.4}{\familydefault}{\mddefault}{\updefault}{\color[rgb]{0,0,0}and store Derivatives}%
}}}}
\put(6076,-2536){\makebox(0,0)[lb]{\smash{{\SetFigFont{12}{14.4}{\familydefault}{\mddefault}{\updefault}{\color[rgb]{0,0,0}Assemble residual (51)}%
}}}}
\put(6301,-2761){\makebox(0,0)[lb]{\smash{{\SetFigFont{12}{14.4}{\familydefault}{\mddefault}{\updefault}{\color[rgb]{0,0,0}and its Jacobian}%
}}}}
\put(8519,-4936){\makebox(0,0)[lb]{\smash{{\SetFigFont{12}{14.4}{\familydefault}{\mddefault}{\updefault}{\color[rgb]{1,0,0}Go to new time step}%
}}}}
\put(10182,-1786){\makebox(0,0)[lb]{\smash{{\SetFigFont{12}{14.4}{\familydefault}{\mddefault}{\updefault}{\color[rgb]{0,0,0}Flag cells and}%
}}}}
\put(5656,-1636){\makebox(0,0)[lb]{\smash{{\SetFigFont{12}{14.4}{\familydefault}{\mddefault}{\updefault}{\color[rgb]{0,0,0}Assemble FEM systems (46) (48)}%
}}}}
\put(6875,-8888){\makebox(0,0)[lb]{\smash{{\SetFigFont{12}{14.4}{\familydefault}{\mddefault}{\updefault}{\color[rgb]{0,0,0}End}%
}}}}
\put(6834,2209){\makebox(0,0)[lb]{\smash{{\SetFigFont{12}{14.4}{\familydefault}{\mddefault}{\updefault}{\color[rgb]{0,0,0}Start}%
}}}}
\put(3246,-2011){\makebox(0,0)[lb]{\smash{{\SetFigFont{12}{14.4}{\familydefault}{\mddefault}{\updefault}{\color[rgb]{0,0,0}find new solution guess}%
}}}}
\put(3335,-1800){\makebox(0,0)[lb]{\smash{{\SetFigFont{12}{14.4}{\familydefault}{\mddefault}{\updefault}{\color[rgb]{0,0,0}Invert Jacobian and}%
}}}}
\put(5614,-4076){\makebox(0,0)[lb]{\smash{{\SetFigFont{12}{14.4}{\familydefault}{\mddefault}{\updefault}{\color[rgb]{0,0,0}no}%
}}}}
\put(8326,-5943){\makebox(0,0)[lb]{\smash{{\SetFigFont{12}{14.4}{\familydefault}{\mddefault}{\updefault}{\color[rgb]{0,0,0}no}%
}}}}
\put(6666,-4395){\makebox(0,0)[lb]{\smash{{\SetFigFont{12}{14.4}{\familydefault}{\mddefault}{\updefault}{\color[rgb]{0,0,0}tolerance?}%
}}}}
\put(6483,-4032){\makebox(0,0)[lb]{\smash{{\SetFigFont{12}{14.4}{\familydefault}{\mddefault}{\updefault}{\color[rgb]{0,0,0}Residual norm}%
}}}}
\put(6928,-4216){\makebox(0,0)[lb]{\smash{{\SetFigFont{12}{14.4}{\familydefault}{\mddefault}{\updefault}{\color[rgb]{0,0,0}$<$}%
}}}}
\put(6643,-7430){\makebox(0,0)[lb]{\smash{{\SetFigFont{12}{14.4}{\familydefault}{\mddefault}{\updefault}{\color[rgb]{0,0,0}Ref cycles}%
}}}}
\put(6984,-7610){\makebox(0,0)[lb]{\smash{{\SetFigFont{12}{14.4}{\familydefault}{\mddefault}{\updefault}{\color[rgb]{0,0,0}=}%
}}}}
\put(7086,-6773){\makebox(0,0)[lb]{\smash{{\SetFigFont{12}{14.4}{\familydefault}{\mddefault}{\updefault}{\color[rgb]{0,0,0}yes}%
}}}}
\put(7097,-4823){\makebox(0,0)[lb]{\smash{{\SetFigFont{12}{14.4}{\familydefault}{\mddefault}{\updefault}{\color[rgb]{0,0,0}yes}%
}}}}
\end{picture}%